*The anthyphairetic reconstruction*
*of the original Pythagorean proof of incommensurability,*
*by means of the restoration of Book II of the* Elements
*to its original Pythagorean form*

*by Stelios Negrepontis and Vassiliki Farmaki*

*To Athanase in friendship and admiration:*
*ἄνθρωπος Θαλῆς* [*the man is a Thales*]
Aristophanes, *Birds.* Line 1009

*Abstract*. Unquestionably the greatest discovery of the Pythagoreans is the existence of incommensurable magnitudes, most probably the incommensurability of the diameter to the side of a square, but there is no agreement among historians of Greek mathematics on their method of proof. In this chapter we present novel arguments not only for an anthyphairetic reconstruction of the original Pythagorean proof of incommensurability, but also in favor of one that employs the Pythagorean Application of Areas in Excess and in fact Geometric Algebra. The main tool for this reconstruction is the restoration of Book II of the *Elements* to its original Pythagorean form. Thus we place the Pythagorean theorem at the position II.4/5 (namely immediately after Proposition II.4 of the *Elements*), the Application of Areas in Defect at II. 5/6, in Excess at II.6/7, the Elegant theorem(s) at (II.9/10 and) II.10/11.a, and the proof of the Pell property of the side and diameter numbers by mathematical induction at II.10/11.b. We then naturally expect that the original anthyphairetic Pythagorean proof of incommensurability is placed, *before* Propositions II.9 & 10 on the side and diameter numbers, and *after* the tools needed for its proof, the Application of Areas II.6/7 (and II.8). hence *precisely* at the position II.8/9.
This reconstruction is strongly confirmed by
(a) the interpretation of the Pythagorean philosophical principles Infinite and Finite, in terms of the Pythagorean incommensurability proof, and
(b) the Pythagorean derivation of the three kinds of angles from these principles, of the acute and obtuse angles from the Infinite, of the right angle from the Finite;
and in turn the reconstruction confirms and explains
(a) the use of Geometric Algebra for its proof,
(b) the use of Mathematical Induction for the proof of Pell's property of the side and diameter numbers by means of Proposition II.10 of the *Elements*, and
(c) the early dating of the Pythagorean incommensurability discovery, in connection with our study of Zeno's arguments and paradoxes,
(d) the presence of Proposition II.8 of the *Elements*, and
(e) the late, in fact after Proposition II.10, appearance of the converse to the Pythagorean theorem, Propositions II.12 and 13.

*Introduction.*

In the present chapter we offer new arguments that help resolve some of the most enduring questions on Pythagorean Mathematics and Philosophy. The greatest mathematical discovery of the Pythagoreans, perhaps the greatest ever, is the discovery of the incommensurability (Section 1). There are ancient proofs of the incommensurability: an elementary arithmetical one, a later anonymous addition to the *Elements* as "Proposition X.117", a hint to this arithmetical proof by Aristotle in *Analytics Prior* 41a which probably refers to a more general proof by Archytas, and, as it is not generally realized, because it is couched in philosophical language, an anthyphairetic proof in Plato's *Meno* 80e-86e, 97a-98b, probably due to Theaetetus, the oldest one preserved. In addition, there are modern reconstructions, a modified arithmetical one by Knorr, and anthyphairetic ones by Chrystal, 1889, Tennenbaum, 1950, Fowler, 1994, that turn out to be consequences of Propositions II.9 or 10, and could well be ancient ones.

There are at least two strong, but still *tentative* arguments, according to which an *arithmetical* proof must be *rejected* as the original Pythagorean incommensurability proof, in favor of an *anthyphairetic* proof, namely in favor of a proof that employs Proposition X.2 of the *Elements,* according to which an *infinite anthyphairesis* of two magnitudes is a criterion of their *incommensurability*.

One of these arguments is the Pythagorean stories about a disrespectful Pythagorean, probably Hippasus, which, according to our interpretation, correlate *the incommensurability with the Infinite*; the second argument is based on the fact that the Pythagoreans had introduced and studied *the side and diameter numbers*, whose anthyphairesis is an initial finite segment of the infinite anthyphairesis of the diameter to the side of a square, the modern "convergents" of the continued fraction of the square root of 2 (Section 2).

The anthyphairetic proof of the incommensurability diameter to the side of a square by Theaetetus, which is the model for *Meno*'s theory of knowledge as recollection, employs the theory of ratios of magnitudes, a theory created by Theaetetus, and thus certainly is not the original Pythagorean proof (Section 3).

At this point a heuristic discussion intends to point out that Propositions II.9 & 10 of the *Elements* are in close connection with the side and diameter numbers and, through them, with the infinite anthyphairesis of the diameter to the side of a square, suggesting that the restoration of Book II of the *Elements* to its original Pythagorean form may lead us, *not only* to arguments, rejecting each of the three modern anthyphairetic reconstructions mentioned above as the original Pythagorean proof of incommensurability. *but, even more significantly* to the original Pythagorean incommensurability proof and to the tools needed for its proof (Section 4).

So, we take our starting point from now on to be Book II of the *Elements* (Section 5).

The first partial restoration of Book II to its original Pythagorean form consists in adding the Pythagorean theorem, as Proposition II.4/5, and with the Pythagorean Application of Areas in Defect, as Proposition II.5/6, and in Excess, as Proposition II.6/7 (Section 6).

The additions of these Propositions on Application of Areas in Book II are related to the controversy of Geometric Algebra and are not universally accepted; thus these additions should be regarded as *tentative*, pending on additional arguments to be provided in Sections 12 and 14 (Section 7).

The second partial restoration of Book II to its original Pythagorean form consists in adding as corollaries to Proposition II.10, according to Proclus' *Commentary to Plato's Republic*, the so called *Elegant theorem* as Proposition II.10/11.a, and the Proposition on *the Pell property of the side and diameter numbers* as Proposition II.10/11.b; it is then natural to add a corresponding

*Subtractive Elegant theorem*, as Proposition II.9/10, as the corresponding corollary to Proposition II.9, a proposition twin to Proposition II.10. Even though the arguments, that Proclus in his *Commentary to Plato's Republic* describes a *bona fide* proof of the Pell property of the side and diameter numbers by mathematical induction based on Proposition II.10, are solid and valid, we still choose to still consider them as tentative, as rests is also controversial, hence this addition should also be regarded as *tentative*, pending on an additional argument to be provided in Section 14 (Section 8).

At this point we study the three anthyphairetic reconstructions, by Chrystal, 1889, adopted by Rademacher and Toeplitz and by van der Waerden, by Tennenbaum, about 1950, and by Fowler, 1994, suggested as the original Pythagorean proof of incommensurability. Fowler's proof is based on the Elegant theorem, while Chrystal's and Tennenbaum proofs are found to be based on the Subtractive Elegant theorem (Knorr, 1998 related Chrystal's proof with Proposition II.9).The close association of these three proofs with Propositions II.9 & 10, both of which are closely related to the side and diameter numbers, which according to Proclus' comments, but also according to their nature as the "convergents" of the continued fraction of square root of 2 in modern language, were definitely *after the discovery and proof of the incommensurability*, lead us to reject them as the original Pythagorean proof (Section 9).

The structure of the Pythagorean Book II of the *Elements* now leads us to expect that the original incommensurability proof must be sought *before* the propositions related to the side and diameter numbers, namely at the position II.8/9 (Section 10).

Our expectation proves fruitful in that we obtain a proof of the incommensurability employing the emblematic Proposition II.4 and its consequences, namely the Pythagorean theorem II.4/5, the "completion of the square" II.6, the Application of Areas in Excess II.6/7, and finally Proposition II.8 (not having any known use till then). The Application of Areas in Excess is used in two basic ways for the incommensurability proof, of which one is typically school-algebraic, in the sense of van der Waerden, and provides the arithmetical value 2 for the second and third anthyphairetic quotients, while the second provides the preservation of the Gnomon, that shows that all the succeeding quotients, infinite in multitude have also the same arithmetical value 2 (Section 11).

At this point we realize that Theaetetus' proof, given in Section 3, employed philosophically in the *Meno*, essentially transforms, with the help of Theaetetus' theory of ratios of magnitudes, the awkward Pythagorean Finitization of the preservation of Gnomons to the Logos Criterion for anthyphairetic periodicity.

The validity of our proof as the true reconstruction of the original one is *tentative* at this point, and must be supported by objective arguments in order to be considered the definitive one. We offer two such objective arguments.

The first is our interpretation of the two philosophic Pythagorean principles, the Infinite and the Finite, on the basis of the description of these principles by Aristotle, Simplicius, and Philoponus: the Infinite is the infinite anthyphairetic division of the diameter to the side of a square, while the Finite/Finitizer is the preservation of the Gnomons, providing full knowledge of an infinite and unknown entity (Section 12).

The certain dating of the discovery is subtle. The association of Hippasus in the Pythagorean stories does suggest an early dating but this is a rather weak evidence. Far stronger is the establishment that Zeno's arguments and paradoxes, in particular Zeno's Fragment B3, are closely connected with Pythagorean incommensurability (Section 13).

The second argument in favor of our reconstruction is our interpretation of the Pythagorean treatment of the three kinds of angles. mentioned in Plato's *Republic* 510c4-5 and explained in detail by Proclus, in terms of the two Pythagorean principles, in which the Infinite, the infinite anthyphairetic division, is the cause of the acute and obtuse angle, and the Finite, the preservation of the Gnomons, the cause of the right angle. This second interpretation confirms not only our reconstruction of the original Incommensurability proof, but also the use of Geometric Algebra (for the proof of the Incommensurability) and of Mathematical Induction (for the proof of the Pell property by the side and diameter numbers by means of Proposition II.10) by the Pythagoreans, and the reason why the converse to the Pythagorean theorem, Propositions II.12 and 13, appears late, and in fact after II.10 (Section 14).

Book II in its restored original Pythagorean form achieves *a unity*, that is lacking in the *Elements*. It is no longer a Book consisting of disparate propositions with no discernible common theme, but it is transformed into a treatise aiming to develop the geometrical/algebraic tools, to use them in order to compute the infinite anthyphairesis of the diameter to the side of a square, and to establish the resulting incommensurability, and finally to examine the side and diameter numbers, the anthyphairetic approximations of the side and diameter numbers. Euclid destroyed this unity by moving the Pythagorean theorem to I.47, its converse to I.48, the Application of Areas to VI.28 & 29, and the proof of incommensurability to problematic X.9.

**Section 1.** *Incommensurability, the great mathematical discovery of the Pythagoreans, and the question on the method of proof employed by the Pythagoreans to prove it*

According to Definition X.1 of the *Elements*, two lines (or more generally magnitudes) a, b are called *commensurable* if there is a common measure, a line c, that measures exactly both lines a and b, namely that for some numbers m and n, we have a=mc and b=nc. If a and b are not commensurable, then they are called *incommensurable*.

The greatest mathematical discovery of the Pythagoreans, arguably the greatest mathematical discovery ever, was their rigorous proof that there are lines a and b, such that a and b are incommensurable; in fact they proved that the diameter and the side of a square are incommensurable. The ancient sources that mention this discovery include Proclus' Synopsis of the History of Greek Mathematics, in his *Commentary to the First Book of Euclid's Elements* 65,15-21, a Synopsis that most likely goes back to Eudemus, and the *Anonymous Scholion In Eucliden* X.1, lines 21-47; and also the sources that recount persistently the story of a Pythagorean, most probably Hippasus of Metapontium, who disrespectfully made public to the "unworthy" non-Pythagoreans the secret discoveries of the Pythagoreans on incommensurability, and was punished by the divinities by drowning/perishing in the sea:

Plutarch, *Nomas* 22,3,4-4,1,
Iamblichus, *On the Pythagorean Life* 246,10-247,1-7; 88,13-89,4, and *On the Common Mathematical Science* 25, 27-39,
Pappus, in his *Comments on Book X of Euclid's Elements, and
Anonymous Scholion in Eucliden* X.1, 70-79
(Our interpretation of the story is outlined in Section 2.4, below).

*Question 1*. The greatest question concerning the mathematical achievements of the Pythagoreans is by what method the Pythagoreans obtained the original proof of incommensurability.

**Section 2.** *The arithmetical proof of incommensurability, hinted by Aristotle and fully described as "Proposition X.117" of the* Elements*, with adherents from Heath, 1921 to Netz, 2022, is tentatively rejected as the original Pythagorean proof in favor of an anthyphairetic proof*

**2.1**. *An elementary arithmetical method of incommensurability hinted in Aristotle's* Analytics Prior *41a and fully preserved in "Proposition X.117", a later anonymous addition to the* Elements

Some ancient sources, including Aristotle, *Analytics Prior* 41a and "Proposition X117" of Euclid's *Elements*, describe arithmetical proofs of the incommensurability of the diameter to the side of a square, and this has led a number of historians of Greek Mathematics to suggest that the original Pythagorean proof of incommensurability is to be found in the arithmetical proof of "Proposition X.117", in fact a later anonymous addition to Euclid's *Elements*.

**2.2**. *Several historians of Greek Mathematics are adherents to the arithmetical original method of proof*

Among the adherents of such an original Pythagorean proof are Heath, 1921, Becker, 1936, Szabo, 1969, Zhmud, 2006, Netz, 2022. It is rather amusing to follow Netz, 2022, p.80-81, the last of the adherents, describing in his own words the proof of X.117:

> If so–already going beyond the scope of Plato's *Meno*–the diagonal of the square is what we would call √2 of the side of the square. It appears that Greeks were aware, perhaps very early on, that with the side of the square taken as the basis of measurement, the diagonal is irrational. How do we know that?
> This now becomes *a little like Archytas's music theory*.
> Assume the opposite–that the side and the diagonal are commensurable.
> This is tantamount to the claim that they are in the ratio of two integer numbers, say, m:n.
> Let us also assume that m:n are the smallest such integers standing in this ratio
> *–the key concern for Archytas's arithmetic!*
> So, m:n are not like 6:8,but like 3:4. They are m:n–the smallest such numbers.
> We also know, by Pythagoras's theorem (always the very same basic tool!), that $2m^2=n^2$.
> Now, because n2 is the double of something (it is the double of $m^2$), it is an even number.There is no way for a square number to be even without its root being even as well (odd times odd is odd, not even). So, n itself is an even number. Indeed, because n is an even number, it follows that $n^2$ is divisible by 4. (An even number is twice something; its square is therefore four times the something squared.) But we have $2m^2 = n^2$ or $m^2=n^2/2$. Now, because $n^2$ is divisible by 4, it follows that $n^2/2$ is divisible by 2 or, in other words, is even.
>
> $m^2=n^2/2$, and $n^2/2$ is even.
>
> We have learned, then, that $m^2$ is even. By the logic applied before, we've also learned that m itself is even. We have therefore established that both m and n are even, contradicting our preliminary hypothesis that m and n are the smallest integers of their ratio. (If they are both even, we just need to divide them both by 2 to get smaller integers with the same ratio. If they are both even, they are not like 3:4, but like 6:8.)
Something is incurably wrong about our initial assumption, and in other words, no two integers can be found to display this ratio. The diagonal of the square is irrational.
I do not know who was the first Greek to observe this. The observation is ingenious but simple; it could well have come about in the very first years of Greek mathematical writing. […]
We only know about this result now because it is tucked in as a kind of appendix at the end of Book X of Euclid's *Elements*.

**2.3**. *Some doubts about whether Aristotle was in possession of the simple X.117 proof of incommensurability*

Van der Waerden, 1954, p.110 believed that the proof of incommensurability to which Aristotle refers in *Analytics Prior* 41a essentially coincides with the proof of Proposition X.117. However, the detailed comments of Alexander Aphrodisieus on Aristotle's passage 260,7-261,28 describe a rather different proof of incommensurability of the diameter to the side of a square, the difference being principally its use of the strong Proposition VII.27

The *proof* of incommensurability of the diameter that Alexander describes in his *Commentary to Aristotle's Analytics Prior* 41a is as follows:
Let $d^2=2s^2$. Suppose that d, s are commensurable.
Then there are numbers m,n, such that d/s=m/n.
By *Propositions VII.20-22*, we may assume that m,n are relatively prime.
By the fundamental *Proposition VII.27*, $m^2$, $n^2$ are also relatively prime.
By *Proposition X.9*, $d^2/s^2=m^2/n^2=2/1$,
Hence $m^2$ is even. Hence $m^2/2$ is even, and since $m^2/2=n^2$, it follows that $n^2$ is even.
Thus $n^2$ is both odd and even, a contradiction.

Alexander's proof is striking, in that in employs the unneeded 'heavy artillery' Proposition VII.27 to prove that $n^2$ is both odd and even, whereas Proposition X.117 reaches this conclusion with entirely elementary means. Thus we must conclude that the elementary proof of X.117 was found *after* Alexander's time, possibly by Theon. *A fortiori*, then, the proof of X.117 was *not* available to Aristotle, and hence it is suggested that Aristotle, with his comments in *Analytics Prior* 41a, had in mind a proof similar to that described by Alexander. This suggests, in turn, that the proof that Aristotle had in mind was a specialization of a more general incommensurability theorem, involving VII.27 or one of its consequences, Proposition VIII.8, 14, 15, closely associated with Archytas. Thus the elementary arithmetical proof of X.117 appears to be a final later product, not available yet in Aristotle's time, and thus not the original Pythagorean proof of incommensurability.

**2.4.** *Our interpretation of the various versions of the Pythagorean story involving the disrespectful Pythagorean correlate Incommensurability with the Infinite, suggesting a Pythagorean anthyphairetic proof of Incommensurability by the use of Propostion X.2*

Pappus, in his *Comments on Book X of Euclid's Elements* and the *Anonymous Scholion In Eucliden* X.**1**, lines 70-79, provides a fascinating explanation of this story by suggesting that the sea in which Hippasus was drowned symbolizes the Infinite, clearly referring to Plato's *Statesman* 273d4-e4 "infinite sea of dissimilarity".

> καὶ τότ' ἤδη θεὸς ὁ κοσμήσας αὐτόν,
> καθορῶν ἐν ἀπορίαις ὄντα,
> κηδόμενος ἵνα μὴ χειμασθεὶς ὑπὸ ταραχῆς
> *διαλυθεὶς εἰς τὸν τῆς ἀνομοιότητος ἄπειρον ὄντα πόντον* δύῃ,
> πάλιν ἔφεδρος αὐτοῦ τῶν πηδαλίων γιγνόμενος,
> τὰ νοσήσαντα καὶ λυθέντα *ἐν τῇ καθ' ἑαυτὸν προτέρᾳ περιόδῳ στρέψας,*
> κοσμεῖ τε καὶ ἐπανορθῶν *ἀθάνατον* αὐτὸν καὶ *ἀγήρων* ἀπεργάζεται. *Statesman* 273d4-e4

> Therefore at that moment God, who brought order in the universe,
> perceiving that it was in dire trouble,
> and fearing that it might *founder in the tempest of confusion*
> *and sink in the infinite sea of dissimilarity,*
> he took again his place as its helmsman,
> *reversed* in the previous period whatever had become unsound and unsettled
> when the world was left to itself,
> set the world in order,
> restored it and made it *immortal and ageless*.

[Plato. *Plato in Twelve Volumes*, Vol. 12 translated by Harold N. Fowler. Cambridge, MA, Harvard University Press; London, William Heinemann Ltd. 1921, with modification by the authors]

It must be remarked that the Infinite [sea] in Plato's *Statesman* passage is menacing with chaos *only* if it is not periodic, while periodicity tames the Infinite by making it knowable. It is precisely this key to universal knowledge that the Pythagoreans revered and wanted to keep secret.

Pythagorean justice, as Aristotle describes it in *Nicomachean Ethics* 1132b21- 23, is reciprocal "anti-peponthos" justice, namely justice by a punishment corresponding and analogous to the affront; since Hippasus' punishment was drowning into the Infinite, we may conclude that the affront, the divulging of *Incommensurability*, was an affront to the *Infinite*, thus, connecting Incommensurability with the Infinite.

Proclus, connects the principle of the Infinite with Incommensurability in his *Commentary to the First Book of Euclid's Elements (=In Eucliden)* 6,19-22, and more specifically in 59,15-61,17 where he correlates incommensurability with an *infinite* decreasing sequence of Gnomons.

This connection of Incommensurability with the Infinite makes perfect sense if we consider that Euclid's only criterion for incommensurability is

> *Proposition X.2.*
> The *infinity* of the anthyphairesis of two lines implies the *incommensurability* of these two lines to each other.

By the *anthyphairesis* of two lines a and b, with a>b, the following process of successive reciprocal division is meant:

$$a = k_0 b + c_1 \text{ with } c_1 < b,$$
$$b = k_1 c_1 + c_2 \text{ with } c_2 < c_1,$$
$$c_1 = k_2 c_2 + c_3 \text{ with } c_3 < c_2,$$
$$c_2 = k_3 c_3 + c_4 \text{ with } c_4 < c_3,$$
$$c_3 = k_4 c_4 + c_5 \text{ with } c_5 < c_4,$$
$$\ldots$$
$$c_{n-1} = k_n c_n + c_{n+1} \text{ with } c_{n+1} < c_n,$$
$$c_n = k_{n+1} c_{n+1} + c_{n+2} \text{ with } c_{n+2} < c_{n+1},$$
$$\ldots$$

where
the lines $c_1, c_2, \ldots, c_n, \ldots$ are the *remainders*, and
the numbers $k_0, k_1, k_2, k_3, k_4, \ldots, k_n, \ldots$ are the *quotients*.
The anthyphairesis is
either *infinite,* if it goes on for ever,
or *finite,* if there is a (final) step of the anthyphairesis in which a remainder divides the immediately previous remainder.
A proof that the anthyphairesis of a to b is *infinite* thus provides a proof, an *anthyphairetic proof,* that the two magnitudes a and b are incommensurable.

**2.5**. *The Pythagorean infinite musical anthyphairesis in Philolaus Fragment 6+ Boethius, De Institutione Musica III.8 implies a musical incommensurability*

The Pythagorean correlation of the *Infinite* with the *Incommensurable* appears implicitly in the fascinating Philolaus *Fragment* 6. The Fragment sets up a "musical" anthyphairesis of musical intervals under the operation * of composition, and provides the first three steps of the musical anthyphairesis of the octave 2/1 to the fifth 3/2:

   octave=fifth*fourth, fourth<fifth,
   fifth=fourth*tone, tone<fourth,
   fourth=tone*tone*diesis, diesis<tone,

while the fourth step

   tone=diesis*diesis*Pythagorean komma, Pythagorean komma<diesis is preserved by Boethius, *De Institutione Musica III.8*.
The form of the remainder, alternating between a power of 2 divided by a power of 3 or conversely, immediately shows that the remainder is never equal to the ratio 1/1, and that the musical anthyphairesis is infinite, implying, in analogy to Proposition X.2, that there is no musical interval measuring (under composition) both the octave and the fifth.

**2.6**. *The Pythagorean association of Incommensurability with the Infinite suggests an anthyphairetic original proof of incommensurability*

In view of Proposition X.2 of the *Elements* and the corresponding musical incommensurability implied by Philolaus' *Fragment* 6, the Pythagorean stories, associating the Incommensurability with the Infinite is a strong indication that the original Pythagorean proof of incommensurability was achieved by an *anthyphairetic* rather than an arithmetical proof, namely by employing Proposition X.2 of the *Elements,* according to which if two magnitudes possess infinite anthyphairesis, then these magnitudes are incommensurable.

**2.7**. *The (anthyphaireses of the) side and diameter numbers are finite initial segments of the infinite anthyphairesis of the diameter to the side of a square,*
*an indication  that the original Pythagorean proof of incommensurability is not arithmetical, but anthyphairetic, and*
*the (fundamental Pell property of the) side and diameter numbers are connected with Proposition II.10 of the* Elements,
*an indication that the original Pythagorean proof of incommensurability is connected with Book II of the* Elements.

Another strong indication against the arithmetical and in favor of the anthyphairetic method of the original Pythagorean proof of incommensurability is the fact that the Pythagoreans defined and studied at depth *the side and diameter numbers*, *after* their discovery of incommensurability, as emphasized by Proclus.

**2.7.1**. *The recursive definition of the side and diameter numbers in ancient sources*

Based on the ancient sources for the side and diameter numbers, which include:
Plato, *Republic* 546b-c,
Theon of Smyrna 43,5-43,15, 44,18-45,8,
Iamblichus, *Commentary to Nicomachus' Arithmetic* 91,21-92,3, 92,23-93,6, and
Proclus, *Commentary to Plato's Republic* 2,24,16-25,13 and 2,27,1-29,4,
we know that the Pythagoreans introduced and studied the double sequence of the *side $p_n$ and diameter $q_n$ numbers* defined *recursively* by

$p_1 = q_1 = 1$,
$p_{n+1} = p_n + q_n$, $q_{n+1} = 2p_n + q_n$ for every natural number n.

However, from the recursive definition given by the ancient sources, we would be at a loss to conceive a plausible way in which the Pythagoreans discovered these numbers. The recursive definition is usable on proving properties of these numbers, but unenlightening on now the Pythagoreans thought of them.

*Question 2.* How did the Pythagoreans come up with the recursive definition of the side and diameter numbers?

**2.7.2**. *Two Proclus' passages, In Republic 2,27,1-11 and In Eucliden 61,12-17, clearly imply that the Pythagorean proof of incommensurability of the diameter to the side of a square **is prior to** the Pythagorean definition of the approximating side and diameter numbers.*

Proclus informs us that the Pythagoreans introduced this double sequence *after* their discovery and proof of incommensurability:

[1] Proclus, *In Eucliden* 61,12-17

|  |  |
|---|---|
|  | ὅπου δὲ τὸ σύνεγγυς ἀγαπῶμεν,<br><br>and when we are content<br>with approximations, |
| οἷον *εὑρόντες ἐν γεωμετρίᾳ* τετράγωνον<br>τετραγώνου διπλάσιον,<br>*ἐν ἀριθμοῖς δὲ οὐκ ἔχοντες*<br><br>as having found in geometry<br>a square double a given square,<br>but not having such a square in numbers |  |
|  | ἑνὸς δέοντός φαμεν ἄλλον ἄλλου<br>διπλάσιον ὑπάρχειν,<br>ὥσπερ τοῦ ἀπὸ τῆς πεντάδος<br>ὁ ἀπὸ τῆς ἑπτάδος<br>διπλάσιος ἑνὸς δέοντος.<br><br>we say that a square number<br>is the double of another square number<br>short by one,<br>like the square of seven,<br>which is one less<br>than the double of the square of five. |

[2] Proclus, *In Republic* 2,27,1-11

|  |  |
|---|---|
| Ὅτι ἐπειδὴ *ἀδύνατον*<br>ῥητὴν εἶναι τὴν διάμετρον<br>τῆς πλευρᾶς οὔσης ῥητῆς<br>(*οὐ γάρ ἐστιν*<br>*τετράγωνος ἀριθμὸς*<br>*τετραγώνου διπλάσιος·*<br>ᾧ καὶ δῆλον ὅτι<br>*ἀσύμμετρά* ἐστιν μεγέθη,<br><br>Since it is *impossible*<br>for the diagonal to be rational<br>when the side is rational<br>(*for there is no* number squared |  |
| that is twice a square<br>– a fact that makes it clear that<br>magnitudes are *incommensurable* –<br>… | |
| | *ἐπενόησαν* οὕτω λέγει<br>*οἱ Πυθαγόρειοι καὶ Πλάτων,*<br>*τῆς πλευρᾶς* οὔσης *ῥητῆς*<br>*τὴν διάμετρον ῥητὴν*<br>*οὐχ ἁπλῶς,*<br>*ἀλλ' ἐν οἷς δύνανται τετραγώνοις,*<br>*τοῦ διπλασίου λόγου,*<br>*ὃν δεῖ τὴν διάμετρον ποιεῖν,*<br>*ἢ μονάδι [δέ] ουσαν*<br>*ἢ μονάδι πλεονάζουσαν·*<br>πλεονάζουσαν μὲν ὡς τοῦ Δ τὸν Θ,<br>δέουσαν δὲ ὡς τοῦ ΚΕ τὸν ΜΘ.<br><br>the Pythagoreans and Plato<br>contrived a way to say that<br>when the side is rational,<br>the diagonal is not rational simpliciter,<br>but rather lies in the squares<br>of the double ratio<br>which it is requisite<br>for the diagonal to make,<br>falling short or exceeding it by a unit –<br>exceeding [by a unit] in the ratio of 4 to 9,<br>but falling short in the ratio of 25 to 49). |

[Proclus, *Commentary on Plato's Republic*, Volume II, Essays 7-15, translated with an introduction and notes by D. Baltzly, J.F. Finamore, G, Miles, Cambridge University Press, 2022]

**2.7.3**. *The fact that the Pythagoreans introduced and studied the side $p_n$ and diameter $q_n$ numbers is a strong argument against the arithmetical and in favor of an anthyphairetic original Pythagorean proof of the incommensurability of the diameter d to the side s of a square, because Anth(d,s)=[1,2,2,2,…] and Anth($q_n$,$p_n$)=[1,2,…,2 (n-1 times)]*

The fact that the Pythagoreans introduced the side and diameter numbers after their discovery of incommensurability as rational approximations of the diameter to the side of a square, as Proclus makes clear (Section 2.7.2), together with the fact that the following Proposition, not mentioned in ancient sources, reveals that this approximation is anthyphairetic. Noting also as Fowler (1982; 1999 p.240-262) shows that practically all approximating/simplifying expressions in Greek Mathematics were anthyphairetic make it certain that the initial Pythagorean definition of the side and diameter numbers was anthyphairetic and thus occurring *after* they had obtained the full knowledge of the infinite anthyphairesis of the diameter to the side of a square.

*Proposition*

> For every natural number n, the pair of side $p_n$ and diameter $q_n$ numbers
> is unique such that:
> Anth($q_n$, $p_n$) = [1,2,2,2,…,2 (n-1 times)] and
> $q_n$, $p_n$ are relatively prime.

*Proof.*
We note that:
Anth($q_1$, $p_1$)=[1], Anth($q_2$, $p_2$)=[1, 2].
We proceed by induction.
Assume that
Anth($q_n$, $p_n$) = [1,2,2,2,…,2 (n-1 times)]
and $q_n$, $p_n$ are relatively prime.
Thus we have
$q_n = p_n + c_1$
$p_n = 2c_1 + c_2$
$c_1 = 2c_2 + c_3$
…
$c_{n-3} = 2c_{n-2} + c_{n-1}$
$c_{n-2} = 2c_{n-1}$
and $c_{n-1} = 1$,
since $q_n$, $p_n$ are assumed to be relatively prime.
From the first two relations, it follows that
$p_n = 2(q_n - p_n) + c_2$
Thus Anth($p_n$, $q_n - p_n$)=[2,2,…(n-1 times)] and $p_n$, $q_n-p_n$ are relatively prime.
Then Anth($p_n + q_n$, $p_n$)= [2,2,…(n times)] and $p_n + q_n$, $p_n$ are relatively prime.
Then Anth($q_n + 2p_n$, $q_n + p_n$)= [1,2,2,…(n times)] and $q_n + 2p_n$, $q_n + p_n$ are relatively prime,
namely
Anth($q_{n+1}$, $p_{n+1}$)= [1,2,2,…(n times)] and $q_{n+1}$, $p_{n+1}$ $q_n$ are relatively prime.

Schematically, the inductive argument from n to n+1 proceeds in the following steps:

| Anth($q_n$,$p_n$)= [1,2,2,…,2 (n-1 times)] | | | | Anth($q_{n+1}$,$p_{n+1}$)= [1,2,2,…,2 (n times)] |
|---|---|---|---|---|
| | | | $2p_n+q_n =$ ($p_n+q_n$)+ $p_n$ | $q_{n+1} =$ $p_{n+1} + p_n$ |
| $q_n =$ $p_n+c_1$ | | $p_n+q_n=$ $2p_n+(q_n-p_n)$ | $p_n+q_n =$ $2 p_n + (q_n-p_n)$ | $p_{n+1} =$ $2 p_n + (q_n-p_n)$ |
| $p_n = 2c_1+c_2$ $c_1 = 2c_2+c_3$ … $c_{n-3} = 2c_{n-2}+c_{n-1}$ $c_{n-2} = 2c_{n-1}$ | $p_n = 2(q_n-p_n)+c_2$ $c_1 = 2c_2+c_3$ … $c_{n-3} = 2c_{n-2}+c_{n-1}$ $c_{n-2} = 2c_{n-1}$ | $p_n = 2(q_n-p_n)+c_2$ $c_1 = 2c_2+c_3$ … $c_{n-3} = 2c_{n-2}+c_{n-1}$ $c_{n-2} = 2c_{n-1}$ | $p_n = 2(q_n-p_n)+c_2$ $c_1 = 2c_2+c_3$ … $c_{n-3} = 2c_{n-2}+c_{n-1}$ $c_{n-2} = 2c_{n-1}$ | $p_n = 2(q_n-p_n)+c_2$ $c_1 = 2c_2+c_3$ … $c_{n-3} = 2c_{n-2}+c_{n-1}$ $c_{n-2} = 2c_{n-1}$ |

Thus, answering *Question 2*, this Proposition, in connection with Proclus' passages considered in Section 2.7.2, reveals that the Pythagoreans introduced the side and diameter numbers, as rational anthyphairetic approximations of the diameter to the side of a square, as the ancients precursors of the convergents of modern continued fractions, *after* they had already obtained full knowledge of the infinite anthyphairesis [1, period(2)] of the diameter to the side of a square, and thus had proved by Proposition X.2 the incommensurability of the diameter to the side of a square.

**2.8.** *Tentative rejection of the arithmetical reconstructions in favor of an anthyphairetic one*

At this point we have reached the tentative conclusion that we must reject the arithmetical reconstruction of the original Pythagorean proof of incommensurability in favor of an anthyhairetic one. We might say that Netz, 2022 is quite right in associating this arithmetical proof with Archytas, but quite wrong in thinking that this is the original Pythagorean proof of incommensurability. Additionally, Knorr's, 1975, p. 27, reconstruction must also be rejected, because it is a variation of the arithmetical one, succeeding to correlate cleverly but rather artificially the incommensurability with the Infinite, but still utterly unable to explain the generation of the side and diameter numbers.

**Section 3.** *The anthyphairetic proof of the diameter to the side of a square by Theaetetus which is the model for* Meno *80e-86c, 97a-98b employs the theory of ratios of magnitudes and is certainly not the original Pythagorean proof.*

It is not generally realized that Plato's *Meno* 80e1-86c3, 97a9-98b6 is based on an anthyphairetic proof of the incommensurability of the diameter to the side of a square, as analyzed and interpreted in Negrepontis, 2024, a proof that must certainly be attributed to Theaetetus, as it makes use of his theory of proportion for magnitudes. We describe this proof in Proposition 3.2 below; the following elementary proposition is needed.

**3.1.** *Proposition*
        If a, b, c, d are line segments, such that ad=bc, then Anth(a,b)=Anth(c,d).

*Proof.* Let Anth $(a, b) = [k_0, k_1,…]$.
We proceed by induction.
Thus $a=k_0b+e_1$, with $e_1<b$.
Then $cb=ad=k_0bd+e_1d$
By *Proposition I.44* of the *Elements*, there is a line segment $f_1$, such that
      $e_1d=bf_1$.
Since $e_1<b$, it is clear that $f_1<d$.
Thus $cb=k_0bd+e_1d= k_0bd+ bf_1$; then
      $c= k_0d+ f_1$ with $f_1<d$.
Thus Anth $(a, b) = [k_0,k_1,…] = [k_0,$ Anth$(b, e_1)]$, while
    Anth $(c, d) = \quad\quad\quad\quad [k_0,$ Anth$(d, f_1)]$, and
at this point we have established for the tetrad $b, e_1, d, f_1$ the same condition
      $bf_1= de_1$

that we had assumed for the original tetrad a, b, c, d, and we are to prove that
     Anth(b, $e_1$)= Anth(d, $f_1$).
We continue as in the first step, and we finish the proof by induction.

**3.2.** *Proposition* (*Pythagorean*). If b is a line and a is a line such that $a^2=2b^2$, then a is the diameter of the square with side b, and Anth(a,b)=[1, period(2)].

*Proof (Theaetetus).*
[Part I]
By the (isosceles) Pythagorean theorem, line a is constructed as the diameter of the square with given side b. Clearly a>b, Since $a^2=2b^2$, it follows that a<2b
(indeed, if we had a≥2b, then $a^2≥4b^2>2b^2=a^2$, contradiction).
Thus b<a<2b. Hence,

     a=b+$c_1$, $c_1$<b.

It follows that b>2$c_1$. 3b>2a
(indeed, if we had b≤2$c_1$=2a-2b, then 3b≤2a, hence $9b^2≤4a^2=8b^2$, contradiction). Hence,

     b=2$c_1$+$c_2$.

b>2(a-b), 3b>2a.

[Part II]
Then $c_1$=a-b, hence $c_2$=b-2$c_1$=b-2(a-b)=3b-2a, hence
b.$c_2$=b(3b-2a)=$3b^2$-2ab, and
$c_1^2$=(a-b)(a-b)=$a^2$+$b^2$-2ab=$3b^2$-2ab.
Hence b.$c_2$=$c_1^2$. Since b>$c_1$, and

     b.$c_2$=$c_1^2$,

clearly $c_2$<$c_1$,
thus the second step is indeed anthyphairetic

     b=2$c_1$+$c_2$, $c_2$<$c_1$.

By b.$c_2$=$c_1^2$ and Proposition 3.1, we have
Anth (b, $c_1$) =Anth ($c_1$, $c_2$)
[namely b/$c_1$=$c_1$/$c_2$. according to the Theaetetus' definition of proportion for magnitudes].
Then Anth ($c_1$, $c_2$) =Anth (b, $c_1$) = [2, Anth ($c_1$, $c_2$)]. Hence,
Anth (a, b) = [1, 2, Anth ($c_1$, $c_2$)] = [1,2, Anth (b, $c_1$) = [1,2, 2, Anth ($c_1$, $c_2$)] =…= [1, period (2)].
.
*Notes*. This proof is not explicitly preserved in ancient sources, but is the product of our interpretation of the *Meno*. Accepting this interpretation, this anthyphairetic proof is the most ancient proof preserved (and not the proof hinted in Aristotle's *Analytics Prior* 41a). The geometer who would conceive of such a proof would normally be Theaetetus, since Theaetetus was the first geometer to discover a theory of proportion of magnitudes based on equal anthyphairesis, the pre-Eudoxean theory mentioned by Aristotle in *Topics* 158b. In the *Topics* passage, there is no explicit attribution, but there is something of a consensus among historians of Mathematics that the theory is due to Theaetetus.

The reconstruction of Theaetetus' theory of ratios of magnitudes has been obtained by Negrepontis and Protopapas, 2024.

Theaetetus' proof provides *the complete knowledge* of the diagonal with respect to the side of a square with one and a half anthyphairetic division steps plus the condition $b \cdot c_2 = c_1^2$, equivalently $b/c_1 = c_1/c_2$, an instance of Theaetetus' Logos Criterion (consisting in the *repetition/recollection* of Logos). Since the anthyphairesis of a to b is infinite, by Proposition X.2, a and b are incommensurable.

The computation depends on the representation of the first quotient $c_1$ as a-b and the second quotient $c_2$ as equal to 3b-2a. These representations use the *side and diameter numbers* $q_n$, $p_n$; thus $c_1 = p_1 a - q_1 b = a-b$, $c_2 = q_2 b - p_2 a = 3b - 2a$.

The anthyphairetic proof of the Pythagorean incommensurability extracted from the *Meno* is the oldest known to us from ancient sources, but it is certainly not the original Pythagorean proof, since it employs Theaetetus' theory of ratios of magnitudes.

**Section 4.** *Heuristic discussion: the restoration of Book II of the* Elements *to its original Pythagorean form is expected to reveal the original Pythagorean proof of incommensurability*

At this point we have tentatively rejected the arithmetical in favor of an anthyphairetic method, as the original Pythagorean proof of incommensurability, thanks to a combination of the accounts by the ancient authors Theon, Iamblichus, Proclus and our knowledge of modern continued fractions, according to which the Pythagorean side and diameter numbers coincide with the modern convergents of the square root of 2, and are thus closely connected with the infinite anthyphairesis of the diameter to the side of a square; in addition we shall see in Section 8 that, as explained by Proclus in his *Commentary to Plato's Republic* 2,27,24-28,4, Proposition II.10 of the *Elements* appears to be the tool for the inductive proof of the fundamental Pell property of the side $p_n$ and diameter $q_n$ numbers, namely the property $q_n^2 = 2p_n^2 + (-1)^n$ for all natural numbers n. Thus, even though there is no reference to incommensurability in Book II, it appears that Book II has some connection with the infinite anthyphairesis, and hence with the incommensurability, of the diameter to the side of a square.

Secondly, there are at least three modern anthyphairetic reconstructions of the proof of the incommensurability of the diameter to the side of a square, each a possible candidate for the original Pythagorean proof, (a) by Chrystal, 1889, adopted by Rademacher & Toeplitz 1930, Toeplitz 1949, and van der Waerden, 1954, (b) by Tennenbaum, about 1950, unpublished, made known by Conway, 2006, Conway and Shipman, 2013, and (c) by Fowler, 1994. These reconstructions could well have been realized in ancient times, because, as we shall see in Section 9, a common, not immediately obvious, thread in these reconstructions, is their reliance *either* on *the Elegant theorem* (if $a^2 = 2b^2$, then $(a+2b)^2 = 2(a+b)^2$), as called by Proclus, in his *Commentary to the Republic* 2,27,24-28,4, an immediate corollary to Proposition II.10 of the *Elements*, for the case of Fowler's 1994 reconstruction,

*or* on the corresponding statement we call *Subtractive Elegant Theorem* (if $a^2 = 2b^2$, then

$(2b-a)^2=2(a-b)^2$), an immediate corollary to its twin Proposition II.9 of the *Elements*, for the case of Chrystal, 1889, or Tennenbaum. Knorr, 1998 had realized that Chrystal's proof is related to Proposition II.9.

Thus again, Propositions II.9 & 10 of the *Elements* appear to be closely related to Pythagorean incommensurability.

On the other hand, several earlier historians of Greek Mathematics have suggested that there are good reasons to believe that several fundamental Pythagorean discoveries, including
the Pythagorean theorem and
the Pythagorean Application of Areas,
that are not part of Book II of the Elements, should nevertheless be included in it.

Thus, in our quest for the original proof of the Pythagorean incommensurability, we will examine, in Sections 5,6,7, and 8, Book II of the *Elements* and its restoration/completion to its original Pythagorean form.

**Section 5**. *Book II of Euclid's* Elements

Book II of Euclid's *Elements*, generally considered Pythagorean in origin, consists of just 14 Propositions.
It is based on the Fifth Postulate,
(a) in contrapositive form Proposition I.29 with consequences Propositions I.30, I.32-34, and the construction of the square of a given side Proposition I.46, and
(b) in a direct form with consequences, the first part of the proof of the squaring of every rectilinear area Propositions I.42-45.

> *Definition II.2* of *Gnomon*
> *Proposition II.1*. If n is a natural number, and $a_1, a_2, \ldots, a_n$, b are lines,
> then $(a_1+a_2+\ldots+a_n) \cdot b = a_1 \cdot b + a_2 \cdot b + \ldots + a_n \cdot b$
> *Proposition II.2*. If a, b, c are lines and a=b+c, then $a^2=ba+ca$,
> *Proposition II.3*. If a, b, c are lines and a=b-c, then $a^2=ba-ca$,
> Proposition II.4. If a, b are lines, then $(a+b)^2=$
> $a^2+$*Gnomon* in the square $(a+b)^2$ *about* the square $a^2 = a^2+b^2+2ab$.
> *Proposition II.5*. If a, x are lines and a/2>x, then $(a/2)^2 - (a/2-x)^2 = x(a-x)$.
> *Proposition II.6*. If a, x are lines, then $(a/2+x)^2=(a/2)^2+x(a+x)$.
> *Proposition II.7*. If a, b are lines, then $(a-b)^2 = a^2+b^2-2ab$.
> *Proposition II.8*. If a, b are lines, then $(a+2b)^2=a^2+4b(a+b)$.
> *Proposition II.9*. If a, b are lines with b<a<2b, then $a^2+(2b-a)^2 = 2b^2 + 2(a-b)^2$.
> *Proposition II.10*. If a, b are lines, then $(a+2b)^2+a^2=2(a+b)^2+2b^2$.
> *Proposition II.11*. If a is a line, construct a line b, such that $a^2=ab+b^2$.
> *Proposition II.12*. If a, b, c are the sides of a triangle, with the sides a and b forming an obtuse angle, then $a^2+b^2<c^2$.
> *Proposition II.13*. If a, b, c are the sides of a triangle, with the sides a and b forming an acute angle, then $a^2+b^2>c^2$.
> *Proposition II.14*. (The second part of the squaring of every rectilinear area) If a, b are lines, construct a line m, such that $ab=m^2$.

Proposition II.8 is an easy consequence of Proposition II.4, and is not used anywhere in the *Elements*.

*Question Q3*. What is the role of Proposition II.8? Where is it used?

Propositions II.12 & 13 constitute the converse of the Pythagorean theorem II.4/5; their proof needs Propositions II.4,4/5, and II.7.

*Question Q4*. Why is the converse of the Pythagorean Theorem proved so late in Book II, in particular after Proposition II.10?

We will attempt to gain essential understanding of the Pythagorean Mathematics and Philosophy by successive steps of *restoration* of Book II in its original, Pythagorean form.

**Section 6.** *Partial restoration of Book II in its original Pythagorean form:*
*Proposition II.4/5=the Pythagorean Theorem,*
*Proposition II.5/6=the Pythagorean Application of Areas in Defect,*
*Proposition II.6/7=the Pythagorean Application of Areas in Excess.*

Several historians of Mathematics, including Heath, 1926 have suggested that three more Propositions must have been originally within Book II

[1] *the Pythagorean theorem as Proposition II.4/5.*

Proclus, in his *Commentary to the First Book of Euclid's Elements* 426,11-12, tells us that the proof of the Pythagorean presented in Proposition I.47 of the *Elements*, relying on the "paradoxical" Propositions I.35-41, is Euclid's proof, and therefore, we may conclude that it is not the original Pythagorean one. It has been suggested a long time ago by Bretschneider, 1870, and Hankel, 1874 that the original proof of the Pythagorean theorem must have been in Book II, before II.9, where it is first used. Based on Plato's *Meno* 84d3-85b7, analyzed in Negrepontis and Farmaki, 2019, and Negrepontis, 2024, we have argued that the original position of the Pythagorean theorem is at II.4/5, receiving a proof in terms of the emblematic Proposition II.4. Thus

*Proposition II.4/5 (Pythagorean theorem).*
If a,b,c are the sides of a triangle, with the sides a and b forming a right angle, then $a^2+b^2=c^2$.

[2] *Pythagorean Application of Areas as Propositions II.5/6 in defect and II.6/7 in excess.*

Several historians of Mathematics, including Tannery 1882, Zeuthen 1886, Heath 1926, van der Waerden 1954, have suggested that the role and purpose of the two Propositions II.5 & 6, variations of the emblematic Proposition II.4 in the form of "completing a square", is to prove the two Application of areas in Defect and in Excess, respectively, that these Applications were originally in Book II, but that Euclid moved them to Book VI, where he proved more general versions (in defect VI.28, in excess VI.29). Thus, we add

*Proposition II.5/6. Pythagorean Application of Areas in Defect.*
For a and m lines, with a/2>m, to construct a line x, such that $x(a-x)=m^2$.

*Proposition II.6/7. Pythagorean Application of Areas in Excess*.
For a and m lines, to construct a line x, such that $x(a+x)=m^2$.

**Section 7**. *The controversy on Pythagorean Geometric Algebra*

Tannery 1882, Zeuthen 1886, van der Waerden, 1954 and other historians of Greek Mathematics proposed the hypothesis that the two Propositions on the Application of Areas were really solutions of quadratic algebraic equations in the guise of geometry, and that Book II should be interpreted as *Pythagorean Geometric Algebra*.
This hypothesis was later disputed by some historians of Greek Mathematics, including Szabo 1969, Unguru 1975, Mueller 1981.
Even though we believe that the arguments already given in support of the thesis on Geometric Algebra, including the original ones, the additional ones by van der Waerden 1976, Weil 1978, Rashed 2009, Blasjo 2016, Negrepontis and Farmaki, 2019, and Negrepontis, Farmaki, Kalisperi, 2022, are strong and valid, we nevertheless intend to provide in the present paper some novel arguments in support of Geometric Algebra, and therefore we will consider the additions of Propositions II.5/6 and II.6/7 at the present time as *tentative and not final*, subject to further evidence.

*Question Q5*. Do Propositions II.5/6, 6/7 really belong in Book II? Do Prepositions II,5/6, 6/7 indicate that the Pythagoreans were practicing Geometric Algebra? And if so, is there a Pythagorean use of Application of Areas II.5/6 or II.6/7?

**Section 8.** *The proof by mathematical induction, using Propositions II.10, that the side $p_n$ and diameter $q_n$ numbers satisfy the Pell property $q_n^2 = 2p_n^2 + (-1)^n$ for every natural number n*

Proclus, in his *In Republic* [=*Commentary to Pato's Republic*] 2,28,17-29,1, states that the Pythagoreans, after defining the side and diameter numbers, (i) first verified the first few cases, (ii) then conjectured, and (iii) finally proved that these numbers are solutions of *the Diophantine Pell equation,* namely that

$q_n^2 = 2p_n^2 + (-1)^n$ for every natural number n.

According to some modern scholars (including Freudenthal, 1953, van der Waerden, 1954), the Pythagoreans went on *to prove* the Pell property by an inductive argument using Proposition II.10 of the *Elements*, but more recently other modern scholars (including Unguru 1991, Acerbi 2000) have challenged the thesis that the Pythagoreans had obtained a proof by induction, or in fact any kind of proof of the Pell property, claiming that the only thing that Proclus' account shows is the verification of the Pell property for the first few cases.
Freudenthal and van der Waerden did not have convincing arguments that Proclus' account amounted to an inductive proof, and they did not have an argument that the Greeks were in possession of the modern Peano principle of Mathematical Induction. (The Wedberg, 1955 and Acerbi, 2000 claim that Plato's *Parmenides* 149a7-c3 rule "number of terms=contacts +1" was proved by an appeal to a *bona fide* principle of Mathematical Induction, has been questioned by Negrepontis, TA, because this rule, if properly understood, applies to a finite recursion only).

In our work (Negrepontis and Farmaki 2019, Chapter 8), we have presented some novel arguments showing that the Pythagoreans had indeed a step-by-step inductive proof, but without the principle of mathematical induction. We base our claim on two remarkable linguistic bonds, in Proclus' detailed account, that, to the best of our knowledge, have not been noted before.

| | The first linguistic bond is between | |
|---|---|---|
| the wording of *the heuristic arithmetical idea* of the Pythagoreans of the (Proclus 2,25,9-12) that the sum of the squares $q_n^2+q_{n+1}^2$, being equal exactly to $2p_n^2+2p_{n+1}^2$, makes both of them *true diameters,* | | |
| | | and the wording by Proclus (2,27,24-28,4) of the statement of *the geometrical Proposition II.10* |
| | | An immediate consequence of *the geometric Proposition II.10* is *the geometric Elegant Theorem* |
| | The second linguistic bond is between | |
| | | the wording of the statement of *the geometric Elegant Theorem* (2,27,13-16), |
| And the wording of the first (2,28,17-21), the second (2,28,22-24) and the third (2,28,27-29,1) inductive step for the proof of *the arithmetical Pell property* | | |

According to the first linguistic bond the heuristic Pythagorean idea about the fact that two successive arithmetical diameters behave like true diameters was the inspiration for the proof of the geometric Proposition II.10, whose immediate consequence is the geometric Elegant theorem, which in effect serves as the geometric model for every step of the proof by induction of the arithmetical Pell property for the side and diameter numbers.

*Proof outlined.*
The arithmetical form of II.10 is an immediate consequence of the geometric.
From the recursive definition of the side and diameter numbers and Proposition II.10 arithmetized, it follows that

$$q_{n+1}^2 + q_n^2 = 2p_{n+1}^2 + 2p_n^2.$$

We now proceed by induction. The proposition holds for n=1. Assume that it holds for n, namely that
$$q_n^2 = 2p_n^2 + (-1)^n.$$
Thus $q_{n+1}^2 + (-1)^n = 2p_{n+1}^2$, namely
$$q_{n+1}^2 = 2p_{n+1}^2 + (-1)^{n+1}.$$

We claim that in its original Pythagorean form Book II included the following Definition and Propositions:

> *The Recursive Definition (at the position) II.9¯* of *the side and diameter numbers*:
> $p_1 = q_1 = 1$,
> $p_{n+1} = p_n + q_n$, $q_{n+1} = 2p_n + q_n$ for every natural number n,
> *Proposition II.10/11,a=*the *Elegant Theorem:* if $a^2=2b^2$, then $(a+2b)^2=2(a+b)^2$, and
> *Proposition II.10/11.b*=the Pell property of the side and diameter numbers:
> $q_n^2 = 2p_n^2 + (-1)^n$ for every natural number n.

Although this is not mentioned in some ancient source, the Pythagoreans would surely have no difficulty in deducing from Proposition II.9 the following

> *Proposition II.9/10 (=*the *Subtractive Elegant theorem):* if $d^2=2s^2$, then $(2s-d)^2=2(d-s)^2$.

Despite the relevance of our arguments, we will still consider at the present stage of our paper our claim that in fact Proclus describes an *inductive proof* of the Pell property by means of Proposition II.10, as *tentative* and not final.

*Question Q6*. Did the Pythagoreans know proof by mathematical induction, and did they prove the Pell property of the side and diameter numbers, or did they simply verify a few initial cases of the Pell property?

**Section 9.** *The Chrystal, Tennenbaum, Fowler anthyphairetic proofs of the incommensurability of the diameter to the side of a square, their reliance on Propositions II.9 & 10, and their rejection as reconstructions of the original Pythagorean proof of incommensurability*

In this Section we will consider the Chrystal, Tennenbaum, Fowler anthyphairetic proofs of the incommensurability of the diameter to the side of a square, and we will show that these, and some variations, are *consequences* of the two Elegant theorems II.9/10 and II.10/11.a, and thus consequences of Propositions II.9 and 10. These proofs are then rejected as reconstructions of the original Pythagorean proof of incommensurability because of their close ancient link with Propositions II.9 and 10, propositions that are connected with the side and diameter numbers, which, as we argued, *follow* the discovery of incommensurability.

**9.1.** *Ancient and modern proofs of*
*the Subtractive Elegant theorem II.9/10;* if $d^2=2s^2$, then $(2s-d)^2=2(d-s)^2$.

### 9.1.1. *Proof from II.9.* Immediate

### 9.1.2. *Proof from the theory of Gnomons of Book II*
Assume $d^2=2s^2$. Then $s<d<2s$, and we form the lines $d-s$, $2s-d$.
By Proposition II.7, $(2s-d)^2=4s^2+d^2-4sd=6s^2-4sd$
By Proposition II.7, $2(d-s)^2=2d^2+2s^2-4ds=6s^2-4ds$.
Hence, $(2s-d)^2=2(d-s)^2$.

### 9.1.3. *Third Proof*

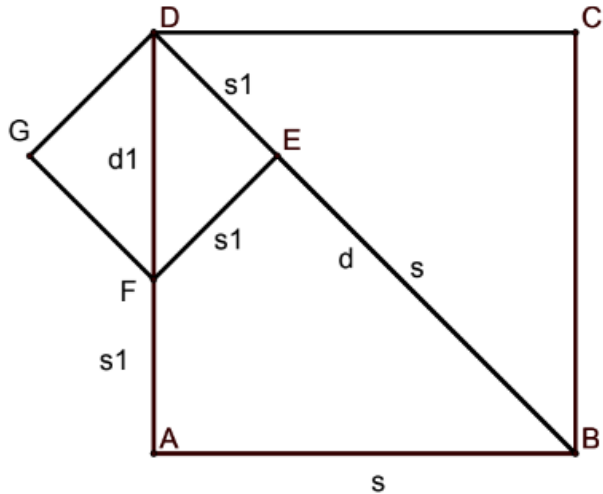

*First step.* Since $d^2=2s^2$, we construct a square ABCD with side $AB=s$ and diameter $BD=d$.
*Second step.* We construct a point E on the diameter such that $BE=s$
Set $ED=s_1$. Then
$d=s+s_1$, hence $s_1=d-s$.
*Third step.* We construct EF perpendicular on diameter DB at point E,
and set F the point of its intersection with AD.
*Claim.* $AF=s_1$
[By Proposition I.34, since DB is a diameter, the angle FDE is equal to half a right angle, and
by construction the angle EDF is a right angle.
By Proposition I.32, the angle DFE is half a right angle.
Hence, by Proposition I.5, $DE=EF=s_1$.]
*Claim.* $AF=FE=s_1$
[The triangles ABF, EBF are right angled, with hypotenuse and one side of the right angle equal]
*Fourth step.* We set $DF=d_1$.
Then DEFG is a square with side $s_1$.
Thus $s_1=d-s$, $d_1=s-s_1=s-(d-s)=2s-d$.
Thus, $s_1=d-s$ and $d_1=2s-d$ are the side and the diameter, respectively, of the square DEFG.
This is the Subtractive Elegant Theorem,

*Note.* Chrystal, 1889, vol. I, p. 270; Rademacher & Toeplitz, 1930, p.16; Toeplitz, 1949, p. 4; van der Waerden, 1954, p. 127; Knorr, 1975, p. 35-36; Knorr, 1998 first realized the close connection of Chrystal's proof with Proposition II.9. We note more specifically that what is needed is just the proof of the Subtractive Elegant theorem.

**9.1.4**. *Fourth proof*. Let s, d be the side and the diameter of a square, so that $d^2=2s^2$. The idea of this proof is to form the great square with side d, and to place inside it the two smaller squares with side s *below left and above right*. Since $d^2=2s^2$, it is clear that the intersection of the two smaller squares must be equal to the sum of the two remaining squares located *above left and below right*.

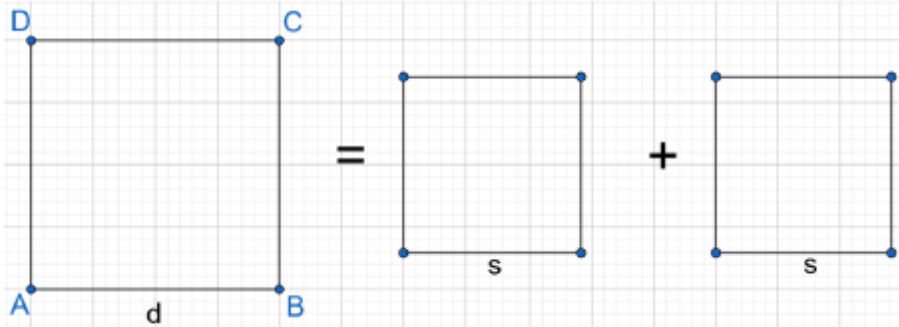

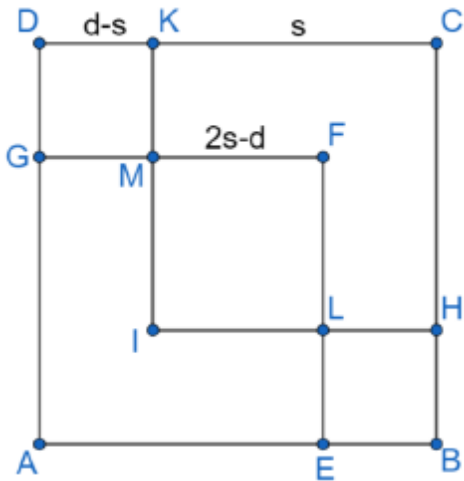

It is easy to see that this equality takes the form $(2s-d)^2=2.(d-s)^2$.
The Subtractive Elegant Theorem follows.

*Note*. Stanley Tennenbaum, about 1950, unpublished, devised this argument for his beautiful proof of the incommensurability of the diameter to the side of a square. This proof became known by Conway, 2006. We realized the close connection of Tennenbaum's proof with Proposition II.9, in fact it is really another proof of the Subtractive Elegant Theorem.

**9.1.5**. *Fifth proof*

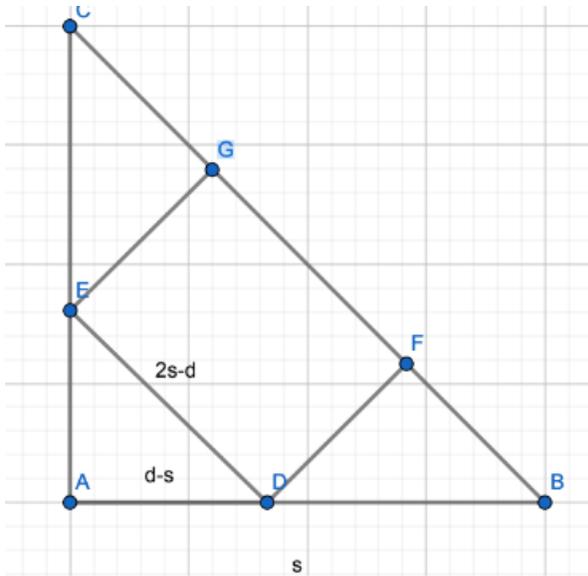

Let ABC be a right triangle with side s=AB=AC and hypotenuse/diameter d=BC.
We set AD=AF=d-s.
We set x=DE, and Claim. x=2s-d.
Indeed,
$x^2 = 2(d-s)^2$ (by the isosceles Pythagorean theorem)
$\quad = 2d^2+2s^2-4ds$ (Proposition II.7)
$\quad = 3d^2-4ds$ (isosceles Pythagorean theorem)
$\quad = d^2+4s^2 -4ds$ (isosceles Pythgorean theorem)
$\quad = (2s-d)^2$ (Proposition II.7).
Thus x=2s-d, and hence
2s-d is equal to the hypotenuse/diameter of the right triangle with sides d-s.

*Note*. We realized that Fowler's argument, given below in 9.3.3, can be turned into still another proof of the Subtractive Elegant theorem.

**9.2**. *The anthyphairetic proof of the incommensurability of the diameter to the side of a square from the Subtractive Elegant Theorem.*

**9.2.1.** *Proposition,* If $d^2=2s^2$, then Anth (d, s) = [1, 2, 2, 2,…],
and thus the diameter and the side of a square are incommensurable to each other.

*Proof* (using *only* the Subtractive Elegant Theorem).

|  | $d^2=2s^2$ |  |  |
|---|---|---|---|
|  |  | Def. $s_1 = d-s$ |  |
|  |  | $d = s+s_1$ |  |
|  |  |  | $d = s+s_1$ |
| Def. $d_1=s-s_1$ |  |  |  |
| $s= d_1+s_1$ |  |  |  |
|  | by the Subtractive Elegant Theorem $d_1^2=2s_1^2$ |  |  |

|   |   |   | Def. $s_2 = d_1-s_1$ |   |
|---|---|---|---|---|
|   |   |   | $d_1 = s_1+s_2$ |   |
|   |   |   |   | $s=2s_1+s_2$ |
|   | Def. $d_2=s_1-s_2$ |   |   |   |
|   | $s_1=d_2+s_2$ |   |   |   |
|   |   | by the Subtractive Elegant Theorem $d_2^2=2s_2^2$ |   |   |
|   |   |   | Def. $s_3 = d_2-s_2$ |   |
|   |   |   | $d_2 = s_2+s_3$ |   |
|   |   |   |   | $s_1=2s_2+s_3$ |
| *ad infinitum* ||||

The general inductive step is as follows:

| Inductive assumption $s_{n-1}=d_n+s_n$ |   |   |   |
|---|---|---|---|
|   | Inductive assumption $d_n^2=2s_n^2$ |   |   |
|   |   | Def. $s_{n+1} =d_n-s_n$ |   |
|   |   | $d_n = s_n+s_{n+1}$ |   |
|   |   |   | $s_{n-1}=2s_n+s_{n+1}$ |
| Def. $d_{n+1}=s_n-s_{n+1}$ |   |   |   |
| $s_n=d_{n+1}+s_{n+1}$ |   |   |   |
|   | by Subtractive Elegant Theorem $d_{n+1}^2=2s_{n+1}^2$ |   |   |

If $s_{n-1}$, $s_n$, $d_n$ have been defined so that $s_{n-1}=d_n+s_n$ and $d_n^2=2s_n^2$,
then, setting $s_{n+1}=d_n-s_n$ and $d_{n+1}=2s_n-d_n$,
we obtain that $s_n=d_{n+1}+s_{n+1}$, and by the *Subtractive Elegant theorem*, that $d_{n+1}^2=2s_{n+1}^2$.
We obtain
Anth $(d,s)=[1,2,2,\ldots]$.
Incommensurability follows from the Proposition X.2.

**9.2.2**. *Note. The arithmetical form of the Subtractive Elegant theorem implies an arithmetical proof of the incommensurability of the diameter to the side of a square*

*Proof.* (S. Tennenbaum's original proof, using the Tennenbaum proof of the Subtractive Elegant Theorem; T.M. Apostol's, 2000 proof using the Chrystal proof of the Subtractive Elegant Theorem). Let s, d be the side and the diameter of a square, so that $d^2=2s^2$, and suppose that d, s are commensurable. Then there are least numbers m, n and a line segment c, such that s=mc, d=nc. We obtain the equation $n^2=2m^2$.
By the arithmetical form of the Subtractive Elegant Theorem, $(2m-n)^2=(2n-m)^2$; and $0 < 2m-n < n$; $0 < n-m < m$; a contradiction, since we have found numbers that express the equality $d^2=2s^2$ and which are strictly smaller than the numbers m, n, respectively.

**9.3.** *Ancient and modern proofs of the Elegant theorem II.10/11.a:* if $d^2=2s^2$, then $(d+2s)^2=2(d+s)^2$

**9.3.1**. *Proof from II.10.* Immediate, as Proclus notes.

**9.3.2.** *Proof from the theory of Gnomons of Book II*
Assume $d^2=2s^2$. We form the lines $d+s$, $d+2s$.
By Proposition II.4, $(d+2s)^2=d^2+4s^2+4sd=6s^2+4sd$
By Proposition II.7, $2(d+s)^2=2d^2+2s^2+4ds=6s^2+4ds$.
Hence, $(d+2s)^2=2(d+s)^2$.

**9.3.3**. *Third Proof* (Fowler, 1994).

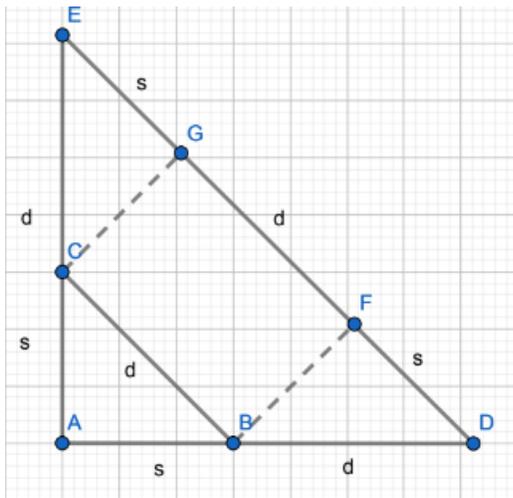

Let ABC be a right-angle isosceles triangle with AB=AC=s, BC=d.
By the isosceles Pythagorean theorem $d^2=2s^2$.
We construct points D,E, such that
A*B*D in same line and BD=d and A*D*E in same line and DE=d,
and we construct the line DE.
Then ADE is an isosceles right triangle with sides s+d.
We construct lines BF and CG perpendicular to the line DE
through the points B, C, respectively..
It is clear that DF=s, EG=s and FG=d. Hence ED=2s+d.
By the isosceles Pythagorean theorem $(d+2s)^2=2(s+d)^2$.

**9.3.4.** *Fourth proof.* (inverse Tennenbaum)

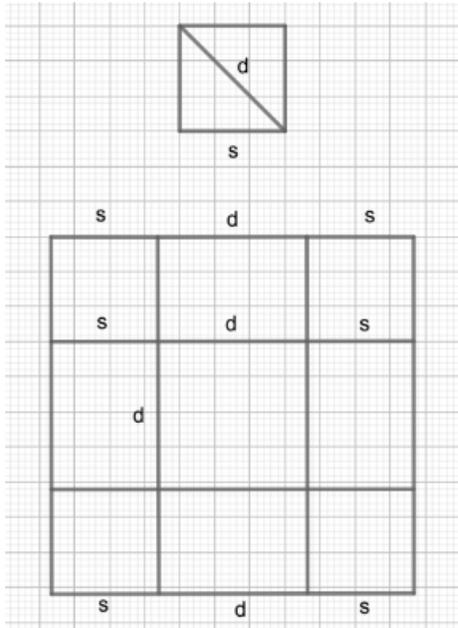

Let a square have side s and diameter d.
We construct the square with side s+d+s; this square consists of
a square $d^2$, four squares $s^2$ and four orthogonal parallelograms sd. Thus
$(s+d+s)^2 =$ (by Proposition II.8)
$d^2+4s^2+4ds =$ (since $d^2=2s^2$)
$2d^2+2s^2+4ds =$
$2(d^2+s^2+2ds) =$ (by Proposition II.4)
$2(s+d)^2$.

**9.3.5.** *Fifth proof* (inverse Chrystal).

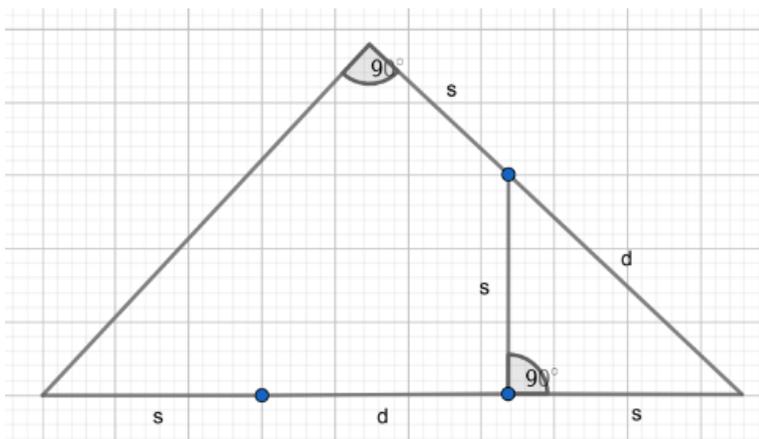

Let a right triangle have right side s and hypotenuse/diameter d.
Extend the side to a line s+d+s,
And extend the hypotenuse d to a line d+s.
We then prove, by reversing the steps in the Chrystal proof (proof 9.3.3), that the line s+d is the right side and the line s+d+s the hypotenuse of a right triangle.

**9.4**. *The anthyphairetic proof of the incommensurability of the diameter to the side of a square from the Elegant Theorem.*

**9.4.1.** *Proposition,* If $d^2=2s^2$, then Anth (d, s) = [1, 2, 2, 2,…],
and thus the diameter and the side of a square are incommensurable to each other.

*Proof* (using *only* the Elegant Theorem).
if $d^2= 2s^2$, namely, if d is the diameter and s is the side of a square,
then $(d+2s)^2=2(d+s)^2$, namely, d+2s is the diameter and d+s is the side of a square;
It is seen immediately that a(a+b)=b(a+2b).
Then by Proposition 3.1, we have
Anth(d,s)=Anth(d+2s, d+s).
On the one hand d>s, and
d=s+(d-s), with d-s<s (since $d^2=2s^2<4s^2$, hence d<2s),
Hence Anth(d,b) = [1, Anth(b, d-b)].

On the other hand, it is clear that
d+2s=(d+s)+s with s<d+s, and d+s=2s+(d-s) with d-s<s.
Hence Anth(d+2s, d+s) = [1,2, Anth(s, d-s)].
Thus [1, Anth(s, d-s)]= [1,2, Anth(s, d-s)].
Hence Anth(s, d-s) = [2, Anth(s, d-s)].
By induction we have Anth(s, d-s) = [2,2,2,…],
and finally Anth(d,s)=[1, 2,2,2,…].
The incommensurability of the diameter to the side follows, by Theorem X.2 of the *Elements,* from the fact that the anthyphairesis of d to s is infinite.

**9.4.2.** *Note*. Although, as we have seen in Section 8, for Proclus the Elegant theorem has the role of a geometrical model for the arithmetical inductive step of the proof of the Pell property of the side and diameter numbers, it is nevertheless true that the Elegant theorem can provide a very simple proof of the computation of the infinite anthyphairesis of the diameter to the side of a square. It is not known if the ancients had ever used the Elegant theorem for this proof. The modern proof was first discovered by D. H. Fowler, 1994.

> 'But, just as Socrates unblocks that impasse [[in the *Meno*]] by conjuring a clever figure out of thin air, so I here draw Figure [1]:
>
> Starting from the small diagonally placed square in the left-hand corner,
> with side s and diagonal d,
> we construct a larger square whose side S is s+d, and check that its diagonal D is 2s+d….

> With one eye on Figure [1], and with the insight that
> the ratios are unaffected by the scale or orientation of the figures involved,
> we take up our problem again, and see that
> the ratio (big side to big diagonal minus side)
> is the same as
> the ratio (little side plus diagonal to little side);
> and we can now evaluate this as
> twice, followed by the ratio (little side to little diagonal minus side)
> which is -scale and orientation aside - what we just started from.
> Hence *the ratio (side to diagonal minus side) is*
> *twice, twice, twice, twice, continuing thus indefinitely,*
> and so the ratio (diagonal to side of a square) is
> once, twice, twice, twice, ....' D. H. Fowler, 1994a, p. 231

**9.4.3**. *Note*. Despite its rudimentary nature, the idea of using the Elegant theorem to calculate the reciprocal of the diameter by the side is not self-evident. Although Proclus does not mention this calculation, it is quite possible that the Pythagoreans were aware of it.
However, this calculation clearly *follows* the discovery of the side-diameter numbers, which demonstrates the crucial relation of the pairs (a, b) and (a+2b, a+b) between "diameter-side" pairs, and this arithmetic recursive relationship is imitated in the geometrical Elegant theorem. So the calculation of the anthyphairesis via the Elegant theorem *follows* the discovery of the diameter-to-side incommensurability, and it must be rejected as the reconstruction of the original Pythagorean calculation of the diameter-to-side anthyphairesis.

***Section 10.*** *The structure of the partially restored Book II of the* Elements *suggests the rejection of the anthyphairetic proofs based on Propositions II.9 and 10, and leads us to the expectation of a reconstruction of the original Pythagorean proof of the incommensurability of the diameter to the side of a square at the position II.8/9.a of Book II*

But now this realization has a fascinating consequence in our quest for the discovery of the original Pythagorean proof of incommensurability. We note that Propositions II.9 and 10 deal with the side and diameter numbers, in fact with the proof of the Pell property of these numbers. But according to common sense, and, in fact, according to the explicit statement by Proclus, the proof of the incommensurability of the diameter to the side of a square must *precede* the introduction, definition and study of the side and diameter numbers, and thus, in the complete restored form of Book II, we expect that the proof of incommensurability of the diameter to the side of a square will *precede* the introduction of side and diameter numbers, namely the incommensurability proof will *precede* Proposition II.9.
We thus expect that the proof of incommensurability of the diameter to the side of a square will be *placed at the position II.8/9*. Schematically:

*Book II partially restored suggests
the position and the original Pythagorean method of proof of
the incommensurability of the diameter to the side of a square*

| II.1,2,3 | | II.9⁻ definition of side & diameter numbers |
|---|---|---|
| II.4 A+b) | | II.9 |
| | | II.9/10 Subtractive Elegant theorem |
| II.4/5 Pythagorean theorem | | II.10 |
| II.5 | | II.10/11.a *Elegant theorem* |
| II.5/6 Application of Areas in Defect | | II.10/11.b Pell property |
| II.6 | | *anthyphairetic definition of side & diameter numbers Proclus, Commentary to Republic 2,27,2-10* II.10/11.c Anth($q_n$,$p_n$)= [1,2,2,…,2 (n-1 times) |
| II.6/7 Application of Areas in Excess | | II.11 Mean and extreme ratio |
| II.7 | | II.12 & II.13 Converse Pythagorean theorem |
| II.8 $(a+2b)^2=a^2+4b(a+b)$ | | II.14 Mean proportional |
| | | |
| | **II.8/9.a** *Original Pythagorean proof of incommensurability before third column, and hopefully with the tools from first column* | |

We are led to the position II.8/9.a for the original position and proof of the Pythagorean incommensurability by looking *backwards* from Propositions II.9 and 10. But if this scheme will have a chance to work only if the tools developed in Book II *before* II.8/9.a, namely the Propositions *up to and including Proposition II.8*, will serve as the tools for this incommensurability proof,

**Section 11.** *The reconstructed original Pythagorean proof of incommensurability at the position II.8/9 of Book II of the Elements, with the use of the Application of Areas in Excess*

*Proposition II.8/9.* If a is the diameter and b is the side of a square then $a^2=2b^2$, Anth $(a, b) = [1,2, 2, 2, …]$, and a, b are incommensurable.

*Proof.* The proof is in four steps

*Step 1.* Since $a^2=2b^2$, we have that $b<a<2b$, we get
$$a=b+c_1, c_1<b.$$
Thus Anth $(a, b) = [1, \text{Anth}(b, c_1)]$.

*Step.2.* We proceed in the anthyphairetic substitution $b+c_1$ for a in the original equality, and get
$$(b+c_1)^2=2b^2.$$
By *Propositions II.4 and II.1,* we find $a^2=((b+c_1)^2=b^2+c_1^2+2bc_1=2b^2$,
$b^2=2bc_1+c_1^2$.

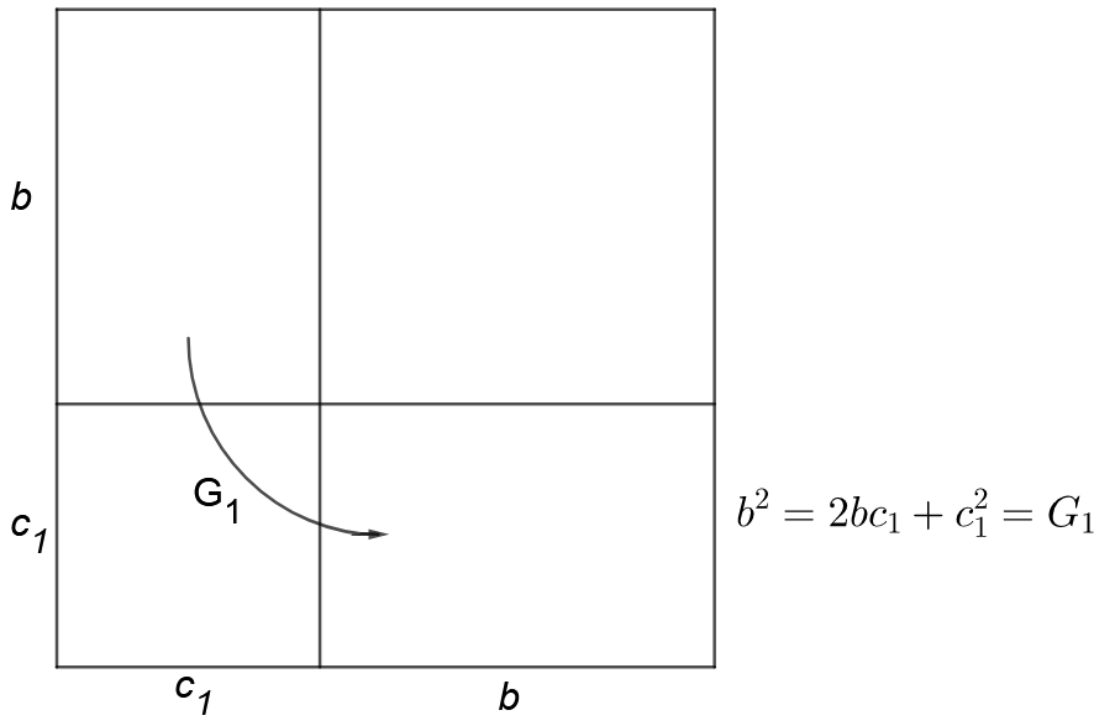

$$b^2 = 2bc_1 + c_1^2 = G_1$$

Then $b(b-2c_1) = c_1^2$.
We set
$$x=b-2c_1,$$
hence $b=2c_1+x$. Then
$$x(2c_1+x) = c_1^2,$$
a Pythagorean Application of Areas in Excess II.6/7
(*in fact the very first use of II.6/7*),
with $a=2c_1$, $m=c_1$.

Thus, the substitution of the first anthyphairetic relation into the original equality results in a Pythagorean application of areas in excess. *Proposition II.6/7* the solution of such an application of areas is the ideal tool for finding, computing the second anthyphairetic quotient.

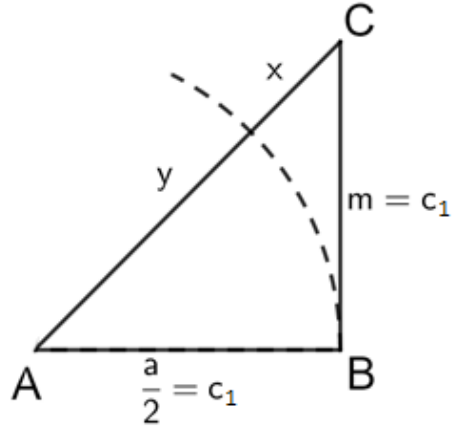

We set **y** for the hypotenuse of the right triangle with right sides $m=c_1$ and $a/2=c_1$.
Then $y^2=c_1^2+c_1^2=2c_1^2$, hence
$\quad\quad y=c_1+c_2$, for some $c_2$, with $c_2<c_1$.
By *Proposition II.6/7*, the solution x is equal to
$\quad\quad x=y-c_1$,
thus $x=y-c_1=(c_1+c_2)-c_1=c_2<c_1$,
hence $b=x+2c_1=2c_1+c_2$, $c_2<c_1$,
$\quad\quad b=2c_1+c_2$, $c_2<c_1$.
Thus, the second anthyphairetic quotient is 2, and we have

$\quad\quad$ Anth $(b, c_1) = [2,$ Anth $(c_1, c_2)]$, and
$\quad\quad$ Anth $(a, b) =[1,$ Anth$(b, c_1)]= [1,2,$ Anth$(c_1,c_2)]$

We note that this use of Application of Areas is purely *"high-school" algebraic*, in the sense meant by van der Waerden. Essentially we solve a quadratic equation, and thus find its integer part.

*Step 3.* Now we proceed in exactly the same way as in the second step.
We substitute b by its anthyphairetic equal $2c_1+c_2$ in the application of areas
$\quad\quad (2c_1+c_2)^2=2(2c_1+c_2)c_1+c_1^2$.

*By Propositions II.8 (first and only use of II.8!),* we find
$c_2^2+4c_1^2+ 4c_1c_2=4c_1^2+2c_1c_2+c_1^2$, then $c_2^2=c_1^2-2c_1c_2=c_1(c_1-2c_2)$, and
$\quad\quad c_2^2=c_1(c_1-2c_2)$.

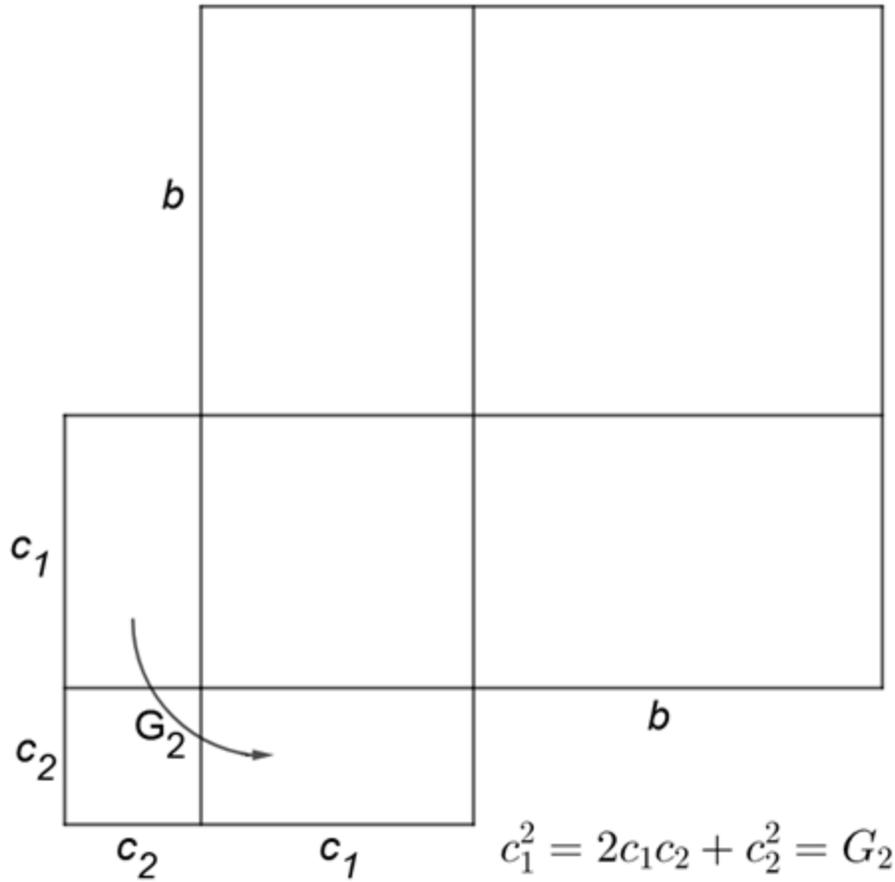

We set
$$x = c_1 - 2c_2,$$
then $c_1 = x + 2c_2$, then
$$x(x + 2c_2) = c_2^2.$$
We observe first that with the anthyphairetic substitution, we have come up again with
*an application of areas in excess II.6/7*,
but secondly that this application of areas
$$x(x + 2c_2) = c_2^2$$
is *the same* as the previous one:
$$x(x + 2c_1) = c_1^2.$$

*Step 4.*
For the completion of the proof we will make use of some propositions making heavy use of the elementary Proposition I.44 of the *Elements* plus a rudimentary induction. We note that Proposition I.44 is in the second part of Book I, which is, as we have already argued, preparatory for Book II, culminating, according to our interpretation with the proof of incommensurability at the exact position II.8/9. Proposition I.44 then finds its use in the proof of Propositions 11.2, 9.4.1, and 11.3.

**11.2.** *Proposition.*
*The preservation of application of areas* (in excess)/ *the preservation of Gnomons.*
   If $Aa^2=Bab+Cb^2$, and $Ac^2=Bcd+Cd^2$, then $ad=bc$.

*Lemma.*
   If $Aa^2=Bab+Cb^2$ and $Aa^2=Bae+Ce^2$, then $b=e$.

*Proof of Lemma.* By hypothesis $Bab+Cb^2=Bae+Ce^2$.
Suppose b unequal to e, say b>e.
Then $Bab>Bae$ and $Cb^2>Ce^2$,
hence $Bab+Cb^2>Bae+Ce^2$, a contradiction.

*Proof of Proposition.*
By *Proposition I.44,* there is a line e such that $ae=bc$.
From the assumption,
$Aa^2c=Babc+Cb^2c$, whence
$Aa^2c=Baae+Cbae$, whence
$Aac=Bae+Cbe$, whence
$Aacc=Baec+Cbec$, whence
$Ac^2a=Bcea+Ce^2a$, whence
$Ac^2=Bce+Ce^2$.
By Lemma, $d=e$.
Hence $ad=bc$.

**11.3.** *Proposition*
   If $Aa^2=Bab+Cb^2$ and $Ac^2=Bcd+Cd^2$, then $Anth(a,b)=Anth(c,d)$.

*Proof.* If $Aa^2=Bab+Cb^2$, and $Ac^2=Bcd+Cd^2$, then $ad=bc$ (Proposition 11.2); and,
if $ad=bc$, then $Anth(a, b)=Anth(c,d)$ (Proposition 3.1).

**11.4.** *The two Pythagorean uses of the Pythagorean Application of Areas in Excess II.6/7 in the original Pythagorean incommensurability proof:*
*(a) Geometric Algebra for Steps 2 and 3, and*
*(b) Preservation of Gnomons for Step 4.*

The *Pythagorean method of the application of areas* has been proved to be crucial in this reconstruction. The application of areas is used in both stages of the proof,
first, in order to construct a solution of the application of areas that is formed at every stage of the division and use this construction/solution in order to calculate the quotient of the anthyphairetic division step, and
secondly, in order to note, when at some stage of the anthyphairesis we come to an application of areas having exactly the same form as in a previous stage, *the preservation of areas/Gnomons,* that the anthyphairesis will be repeated from there on *ad infinitum*.

The *first use,* in our opinion, resolves the long conflict between the supporters and opponents of *Geometric Algebra* in favor of the supporters.

The Pythagoreans/Hippasus wanted *to compute the anthyphairesis* of the diameter to the side of a square so as to prove it is infinite.
They needed to compute *the anthyphairetic quotient* for each one of the infinite steps.
They succeeded in doing this by solving the quadratic equation $x^2=2xy+y^2$. In modern terms setting $\omega=x/y$, $\omega^2-2\omega-1=0$, we find $\omega=(2+\sqrt{8})/2$, hence the integer part of $\omega$ is equal to 2, and this is precisely the anthyphairetic quotient. The Pythagoreans are thus doing algebra precisely in the sense of van der Waerden, 1976, p.199:

> When I speak of Babylonian or Greek or Arab algebra, I mean algebra in the sense of AL-KHWARIZMI, or in the sense of CARDANO'S "Ars magna", or in the sense of our school algebra. Algebra, then, is: the art of handling algebraic expressions like $(a+b)^2$ and of solving equations like $x^2+a\,x=b$.

High School Algebra, indeed:

> Starting with the Application of Areas in excess $b^2=2bc_1+c_1^2$,
> transforming it to modern notation $(b/c_1)^2=2(b/c_1)+1$, and setting $b/c_1=\omega$,
> we obtain the quadratic algebraic equation $\omega^2-2\omega-1=0$.
> Solving it we find
> $\omega=(2+\sqrt{2^2+4.1.1})/2=(2+\sqrt{8})/2=1+\sqrt{2}$, thus $2<\omega<3$.
> Reverting back to ancient notation
> $2/1< b/c_1=\omega<3/1$, namely $2c_1<b<3c_1$, namely $b=2c_1+c_2$, $c_2<c_1$.
> This is precisely the second anthyphairetic step.

Our reconstruction of the original Pythagorean proof of the incommensurability of the diameter to the side of a square with the use of the Pythagorean Application of Areas in Excess II.6/7, a use that is clearly algebraic in the van der Waerden sense of "high school algebra" provides the only use of the method of Application of Areas by the Pythagoreans and, in our opinion, settles the question of *Geometric Algebra*, deciding in favor of the thesis of those historians of Mathematics, such as Tannery, Zeuthen, van der Waerden, who originally proposed it and had the right intuition, but, as far as we know, never explained the use that the Pythagoreans had for the method they discovered.

The second Pythagorean use of the Pythagorean Application of Areas is related to the preservation of the Gnomons, namely that the preservation of the Gnomons, the preservation of the Application of Areas in excess implies the preservation of the anthyphairetic quotient, and thus constitutes a finitization of the infinite anthyphairesis.

**Section 12**. *The anthyphairetic nature of the two Pythagorean philosophic principles Infinite and Finite at the position II.8/9.b provides a confirmation of our reconstruction of the original Pythagorean proof of incommensurability in Section 11*

We will briefly examine the Pythagorean philosophy as presented by Philolaus, Plato, and Aristotle. The most impressive element of the Pythagorean philosophy is that it posits as principles two abstract concepts, the Infinite and the Finite, concepts which derive, according to Aristotle's categorical assurance, from their dealing with "lessons" (Mathematics). We should thus be able to trace, within the Pythagorean mathematical discoveries, the origin of these two principles.

*Question Q7.* What is the meaning of the two Pythagorean principles Infinite and Finite?

Indeed, we will demonstrate that the Pythagorean philosophy is a faithful rendering of the fundamental geometric reconstruction of the Pythagorean anthyphairetic proof of the incommensurability of the diameter to the side of a square, presented in Section 11. Specifically, the two Pythagorean principles Infinite/Apeiron and Finite/Peras correspond to the infinite anthyphairesis of the diameter by the side of a square, and its Finitization by the fact that at each stage of the infinite process the shape and form are preserved, with the subtraction of a Gnomon, in the sense that the anthyphairetic dyad of that stage satisfies the same application of areas in excess.

**12.1**. *According to Aristotle, the two philosophic Pythagorean principles, Infinite and Finite, have mathematical origin*

[1] Aristotle is clear: *the Pythagorean philosophic principles were the principles that were formed by their pre-occupation with Mathematics.*

> Among these thinkers and before them,
> the so-called Pythagoreans were the first to latch on to mathematics.
> They both advanced these inquiries and, having been brought up in mathematics,
> thought that its starting-points were the starting-points of all beings.
> Aristotle, *Metaphysics* 985b23-26

[2] The two Pythagorean principles are *the Finite and the Infinite*,

> The Pythagoreans have spoken of the principles as two … *the Finite* and *the Infinite*
> Aristotle, *Metaphysics*   987a13-16

More carefully Philolaus in his Fragment 1 refers to the two Pythagorean principles as the Infinite and the *Finitizing (περαινόντων).*
Plato in the *Philebus* is crediting some almost mythical Prometheus, who is certainly Pythagoras.

> *Socrates* A gift of gods to men, as I believe, was tossed down from some divine source *through the agency of a Prometheus* together with a gleaming fire; and the ancients, who were better than we and lived nearer the gods, handed down the tradition that all the things which are ever said to exist are sprung from one and many and have inherent in them *the finite and the infinite*.
> Plato, *Philebus* 16c5-10

We then expect to find in the Pythagorean Mathematics the source of these two fundamental Pythagorean philosophical principles.

**12.2.** *The Pythagoreans gave the symbolic names Even/Artion for the principle of the Infinite/Apeiron, and Odd/Peritton for the principle of the Finite/Peras*

The Pythagoreans liked to express themselves symbolically; thus the Finite was called Odd and the Infinite Even. The linguistic similarity of Finite (Peras) with Odd (Peritton) and of Infinite (Apeiron) with Even (Artion) is lost in English.

> the odd and the even,
> and that of these the odd is the Finite and the even the Infinite
> Aristotle, *Metaphysics* 986a18-19

Our task is to obtain from the symbolic descriptions the real content of the two Pythagorean principles.

**12.3.** *Aristotle and Simplicius clarify the symbolic description of the Pythagorean principle of the Finite as Odd: an odd number of units in the figure of an arithmetical Gnomon surrounding a square and thus forming a larger square, thus preserving the square figure; the preservation of the square arithmetical figure during this process* ad infinitum *is precisely the symbolic finitization*

> But there are things that increase *without altering,*
> as the square
> by surrounding it with a *gnomon*
> is increased *but is not thereby altered;*
> Aristotle, *Categories* 15a30-31

The Odd number possesses a finitizing power, in the sense that an Odd number in the form of an arithmetical Gnomon, by surrounding and comprehending a square, increases the square, but preserves the form of the square. Simplicius explains in detail, in his *Commentary to Aristotle's Physics* 456,16-457,8, and we paraphrase:

Why the principle of the *Finite* is symbolized by *the Odd*?
Every *Odd* number 2n+1 can be represented as
an arithmetical *Gnomon* $G_n$ consisting of
n horizontal unit +1 corner unit +n vertical units.
And what is the Infinite we may ask that is finitized by this Gnomon, and in what sense?
The Infinite is the infinite collection of all squares
The connection of "Gnomon" and "Peras-Peritton" is understood, always on a symbolic level, if we note that the arithmetical Gnomon $G_n$, with which we go from a small arithmetical square $T_n$ to the immediately larger square $T_{n+1}$, is always an odd number and has the shape of the Gnomon.
That is, for each arithmetical square
$T_1=1\times1$, $T_n= n\times n$,
the Gnomon/Odd
$G_n =T_{n+1}-T_n =2n+1$,
"comprehends" (εναπολαμβάνει) the square, by "surrounding" (περιθέσεως περί) it,
and thus generates the immediately greatest square
$T_{n+1}=(n+1)\times(n+1)$

in every step of an infinite process.
Thus, the symbolic name Odd for the principle of the Finite is explained by the fact that it is a Gnomon, and the finitizing action of the Gnomon is clearly its comprehension and surrounding of a square, by which it has the power to preserve the shape of the square throughout infinite process. The new square is greater than the old square, but it has not changed in the sense that it retains the same shape and form.

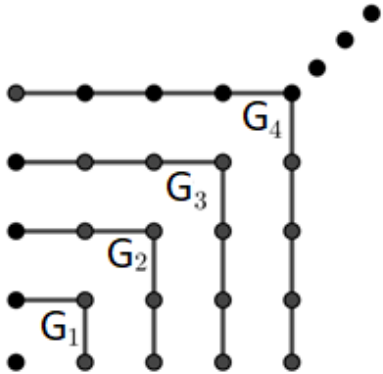

The Odd/Gnomon $G_n$ surrounding the square $T_n$ increases to the square $T_{n+1}$ but does not alter its shape, since after the circumscription of the square with the Odd/Gnomon, the resulting figure is again a square, exactly the same in form and shape as the original smaller square.

**12.4.** *The symbolic description of the Infinite as Even implies that the real Pythagorean principle of the Infinite is a principle of binary division* ad infinitum*,
and in fact, as Simplicius and John Philoponus make clear,
a principle of binary division* ad infinitum *of magnitudes into equal or unequal parts*

> ἅμα γαρ [οι Πυθαγόρειοι] ουσίαν *ποιούσι το άπειρον και μερίζουσιν.*
> With the same breath they [the Pythagoreans] treat the infinite as substance, and divide it into parts.

According to Aristotle, *Physics* 204a33-34, the Pythagorean principle of the Infinite is in fact a Principle of division *ad infinitum*.
Simplicius, *Commentary to Aristotle's Physics* 455,20-23 writes:

> These [the Pythagoreans] said that the infinite was the even number,
> because, according to their expositors, everything even can be divided into equals,
> and what is divided into equal parts is infinite through bisection.
> For division into equal halves proceeds *ad infinitum*

But shortly after Simplicius specifies that the principle of the Infinite is a principle of binary division of magnitudes, not of numbers, and in fact in unequal parts.
> So in this way the expositors ascribe the Infinite to the Even
> through the division into equal parts,
> and it is clear that they take the division *ad infinitum*
> not as being of numbers but *of magnitudes* (οὐκ ἐπ' ἀριθμῶν ἀλλ' ἐπὶ μεγεθῶν).

> For numbers, including the even, do not all divide into equals for long,
> and those that are divided reach the one and halt the division,
> while in the case of *magnitudes* what prevents *the remainder* (τὸ καταλειπόμενον, 455,29)
> from being divided *ad infinitum*, *even if the division is not into equal parts*? 455,24-30

Finally, Simplicius states that in fact Aristotle does not have in mind division into equal parts, but rather into unequal parts.

> In general *Aristotle does not seem to ascribe being Infinite
> to the division into equal parts.*455,33-35

*The similar comments by John Philoponus on the Infinite*

A similar comment to Aristotle's *Physics* 202b30-203a16 is found by John Philoponus, *Commentary on Aristotle's Physics*, 386.14-395.7, especially in 388.24-389.20.
According to this Commentary, as the Pythagoreans presented their teachings
"symbolically" ("συμβολικώς"), "from the numbers", "εκ των αριθμών",
so also the principle "the Infinite" was described "symbolically", as "Even",
something that already becomes clear from what Aristotle says.
For, Even "pertains to division" ("προς της διαιρέσεως εστί") and
> "it is the cause of infinite division", as the dichotomy
> "it is the cause of the division of the magnitudes *ad infinitum*"

and the Infinite is
> "the cause*"(η αιτία)* "of disunity, indeterminacy, division *ad infinitum,* and multitude"
> (της διαστάσεως και αοριστίας και *της επ' άπειρον τομής και πλήθους),*
> "of the magnitudes and of all the entities in nature" (*τοις μεγέθεσι* και πάσι τοις φυσικοίς).

This brings the Pythagorean principle of the Infinite close to a philosophical form of infinite anthyphairesis.

**12.5.** *The symbolic arithmetic Gnomons refer to the real geometric Gnomon of Book II of the Elememts, and thus the two Pythagorean principles arose from the Mathematics in Book II*

Simplicius, *Commentary to Aristotle's Physics*
> The Pythagoreans called *the odd numbers gnomons*
> because when placed around squares they preserve the same shape,
> *like the gnomons in geometry.* For they call the two complements together with one of the parallelograms around the same diagonal *a gnomon*, which, when added to the other of the parallelograms around the same diagonal *makes the whole similar to that to which it was added*.
> So *the odd numbers* are also called *gnomons*
> because, when added around the already existing squares,
> they always preserve the square shape. 457,1-8

The arithmetical Gnomons that symbolize the Pythagorean principle of the Infinite is correlated with the geometric Gnomons of Book II of the *Elements*

**12.6**. *The geometric Gnomons in the real principle of the Infinite must be, not increasing as in the symbolic description, but decreasing, as hinted by Aristotle and specified explicitly by Proclus*

Since the *actual* Pythagorean Infinite refers to binary infinite *division* of magnitudes
(and not to the *symbolic* incremental infinite sequence of numbers),
what would the *actual* Pythagorean Finite/Peras be?
The finitization of the symbolic incremental infinite of the sequence of numbers
is carried out with Odd arithmetical Gnomons which preserve the shape unchanged throughout the infinite process.
Since the Pythagorean principle of the Infinite is
*symbolically* about arithmetical Gnomons increasing by addition *ad infinitum,*
but actually about *geometric Gnomons* decreasing by division *ad infinitum;*
and the principle of the Finite acts upon the Infinite by preserving, not the arithmetic, *but the geometric* Gnomons.

> *it follows that the Gnomons that preserve the form in the real Pythagorean Infinite must be not arithmetical, but geometric, like the Gnomons of Book II of the Elements;*

and, further, since the Infinite is a dyadic division *ad infinitum* into decreasing magnitudes,
it also follows that

> *the geometrical Gnomons that preserve the form in the actual Pythagorean Infinite must be not increasing but decreasing*.

Decreasing Gnomons preserving the form are indeed foreseen in a passage in Aristotle's *Categories*:

> Thus alteration would be distinct from the other changes.
> For if it were the same, a thing altering would, as such, have to be
> *increasing too or decreasing,* or one of the other changes would have to follow;
> but this is not necessary.
> Equally, a thing increasing-or undergoing some other change-would have to be altering.
> But there are things that *increase without altering,*
> *as a square is increased by the addition of a gnomon but is not thereby altered*;
> similarly, too, with other such cases.
> Hence the changes are distinct from one another.
> Aristotle, *Categories* 15,25-33

Proclus' passage, in his *Commentary to the First Book of Euclid's Elements* 59,15-61,17, refers explicitly to the presence of an infinite sequence of descending Gnomons, a feature characteristic of incommensurability:

> The statement that every ratio is commensurable belongs to arithmetic only
> and not to geometry, for geometry contains incommensurable ratios.
> Likewise, the principle that there is a least for the division of square gnomons is peculiar
> to arithmetic: in geometry there is no least gnomon. 60, 7-12

Here Proclus introduces a criterion of incommensurability,
not simply by means of anthyphairetic infinity,
but, more specifically, by the presence of an infinite decreasing sequence of Gnomons.

**12.7.** *Anthyphairetic interpretation of*
*the two Pythagorean principles Infinite and Finite/Finitizing*

According to Sections 6.1-6, behind the Pythagorean symbolism,
*the real Pythagorean principle of the Infinite* is the decrease by binary division *ad infinitum* of magnitudes, in fact of squares and surrounding Gnomons, and
*the real Pythagorean principle of the Finite* is the preservation, despite the *ad infinitum* decrease of magnitudes, of the same form and shape of the squares by the surrounding Gnomons, and in this sense the finitization of the Infinite.

At this point we cannot fail to recognize that a mixture of Infinity and Finite, with characteristics in striking agreement with what is required, appears in one and only one case in Pythagorean geometry: in our reconstruction of the proof of diameter-to-side square incommensurability, in Section 11, where
the Infinite is the infinite process of anthyphairesis of the diameter by the side of a square, and
the Finite is the preservation of the shape and form of the application of areas in excess/Gnomon, at each step.

| *Preservation of Gnomons/Finite* | *Anthyphairetic Division/Infinite* |
|---|---|
| $a^2=2b^2$ | |
| | $a=b+c_1$ |
| $b^2=c_1.(2b+c_1)=G_1$ | |
| | $b=2c_1+c_2$ |
| $c_1^2=c_2.(2c_1+c_2)=G_2$ | |
| | $c_1=2c_2+c_3$ |
| $c_2^2=c_3.(2c_2+c_3)=G_3$ | |
| | $c_2=2c_3+c_4$ |
| … | |
| $c_{2n-2}^2= G_{2n-1}$ | |
| | $c_{2n-2}=2c_{2n-1}+c_{2n}$ |
| $c_{2n-1}^2 =G_{2n}$ | |
| | $c_{2n-1}=2c_{2n}+c_{2n+1}$ |
| … | |

The Finitization of the Infinite in the *anthyphairesis* of the diameter to the side of a square is analogous to the symbolic-arithmetic one, in the following sense:

In the symbolic description with the addition of the Gnomon, the shape and form of the square is increased but preserved indefinitely,
in the infinite multiplication of a diameter by a side with the removal of the Gnomon, the shape and form of this equation of the parabola of sections in hyperbola is preserved indefinitely.
We conclude that the Pythagoreans did receive, as Aristotle expressly informs us, the two principles Apeiron and Peras from their engagement with the lessons, and specifically from their most fundamental discovery, that is, the anthyphairetic proof of the incommensurability of the diameter and side of a square.

**12.8.** *According to Philolaus, Simplicius, Hesychius, Eustratius, Anonymon Scholion to Elements, the Gnomon possesses the power to provide knowledge*

An additional confirmation of our interpretation is provided by the numerous ancient comments on the power of the Gnomon to provide knowledge.

According to Philolaus' (possibly not genuine) Fragment 11
the gnomon has cognitive power.

>    πάντα γνωστὰ καὶ ποτάγορα ἀλλάλοις κατὰ γνώμονος φύσιν ἀπεργάζεται

In the passage of the geometric number in *Republic* 546b3-c7, which has a close relationship with Philolaus' *Fragment 11*, the position of the guide is taken by Logos, which also has cognitive power:
> πάντα προσήγορα καὶ ῥητὰ πρὸς ἄλληλα ἀπέφηναν· *Republic* 546b7-c1
> render all things conversable and rational with one another,

Simplicius, *Commentaries on Categories* 430.5-431.5, states that the geometric gnomon, consisting, according to Definition II.2 of Euclid's Elements, of the two complements with one of the parallelograms

> γνωματεύει και κρίνει
> makes known and judges the remaining rectangle.

Hesychius, *Lexicon* Γ 748 and Eustratios, *In Nicomachean Ethics* 23,15-30
associate Gnomon with knowledge.
*Anonymous Commentary II.11 on Euclid's Elements*, which is a comment on Definition II.2 of the geometric Gnomon, justifies the name of the Gnomon by saying that

> απ' αυτού…το όλον γνωρίζεται
> by Gnomon...the whole becomes known

just as the astronomical Gnomon has as his sole mission

> το τας ενεστώτας ώρας ποιείν γνωρίμους
> The Gnomon is that which makes known the present hours

The cognitive nature and power of the Pythagorean principle is in complete agreement with the role of the geometric principle of the preservation of Gnomons in the reconstruction of incommensurability given in Section 11, where the acquisition of knowledge of the infinite anthyphairetic expansion of the diameter to the side of a square is clearly achieved because of the finitizing power of Gnomon.

So the cognitive nature and power of the gnomon in Pythagorean philosophy is the transference of the corresponding cognitive nature and power of the gnomon in Pythagorean geometry, and especially in the Pythagorean proof of the incommensurability of diameter to side of a square. This association of the proof method, i.e. the anthyphairetic method, with the acquisition of complete knowledge of diameter to side is an additional substantial argument in favor of the anthyphairetic reconstruction, as the traditional arithmetic incommensurability method of proof is far from providing this complete knowledge.

***Section* 13.** *The two philosophical principles Infinite and Finite reappear in Zeno's Fragment B3, with substantially the same anthyphairetic meaning, establishing an early dating of the Pythagorean proof of Incommensurability and a further confirmation of their method of proof*

*Question Q8. Dating of the discovery of the Pythagorean incommensurability*

Democritus, born about 470 BC, wrote 'on irrational lines and solids;' and as Hardy & Wright, 1938, p.42, have observed

> 'it is difficult to resist the conclusion that the irrationality of the diameter to the side of a square was discovered before Democritus' time.'

But the relation of the Pythagorean discovery of incommensurability with Zeno of Elea is even stronger and older. We will be content here to consider

> Καὶ τί δεῖ πολλὰ λέγειν, ὅτε καὶ ἐν αὐτῷ φέρεται τῷ τοῦ Ζήνωνος συγγράμματι;
> πάλιν γὰρ δεικνύς, ὅτι
> εἰ πολλά ἐστι, τὰ αὐτὰ *πεπερασμένα* ἐστὶ καὶ *ἄπειρα*,
> γράφει ταῦτα κατὰ λέξιν ὁ Ζήνων·

*Zeno's Fragment B3*
> "εἰ πολλά ἐστιν,
> ἀνάγκη *τοσαῦτα* εἶναι *ὅσα* ἐστὶ καὶ οὔτε πλείονα αὐτῶν οὔτε ἐλάττονα.
> εἰ δὲ *τοσαῦτά* ἐστιν *ὅσα* ἐστί, πεπερασμένα ἂν εἴη.
>
> εἰ πολλά ἐστιν, *ἄπειρα* τὰ ὄντα ἐστίν.
> ἀεὶ γὰρ *ἕτερα μεταξὺ* τῶν ὄντων ἐστί, καὶ πάλιν ἐκείνων *ἕτερα μεταξύ*.
> καὶ οὕτως *ἄπειρα* τὰ ὄντα ἐστί."
> Simplicius, *Commentary to Aristotle's Physics* 140,27-33
>
> And why should I say any more, for it also exists in the treatise of Zeno?
> For again, showing that if there are many,
> the same things will be *Finite and Infinite*,
> Zeno writes thus verbatim:

'If there are many,
necessarily they are as many as they are, and neither more of them nor fewer.
But if they are as many as they are,
they would be *Finite*.

If they are many,
existing things will be *Infinite*.
For there are always *other* things *generated by* existing things,
and again *other* things *generated by* them.
And in this way existing things will be *infinite*.'

According to the analysis of Zeno's arguments and paradoxes,
undertaken by Negrepontis, 2023,
(a) Zeno's arguments and paradoxes, including Zeno's *Fragment* B3, aim to establish,
by a contradiction,
*not* the impossibility of motion or the impossibility of multiplicity, as is generally believed,
*but* that the true [Plato's intelligible] Beings and the Beings of Opinion [Plato's sensibles] are
separate and different entities.
(b) the true Beings satisfy the almost contradictory property *Infinite and Finite*, namely of
having at the same time an *infinite multitude* of parts, by the principle of Otherness, and
being *One*, in the self-similar sense, according to which each of these infinitely many parts has
*as many* (*tosa*) parts *as (hosa) of* any other part, namely *the same finite number* of parts, while
the sensible entities do not satisfy the *Infinite and Finite* property;
(c) The contradiction results from the assumption that the true Beings coincide with the
sensibles, called the many.
(d) The One of the second hypothesis in the Platonic dialogue *Parmenides*,
(d1) is the philosophical analogue of *periodic anthyphairesis*,
--of *infinite anthyphairesis,* by means of the generation ad infinitum of parts, by the principle of
Otherness *(Parmenides* 143b1-8),
exactly as with the generation of the Infinite in Zeno's Fragment B3, and
--*of circularity* and *periodicity of the anthyphairesis*, resulting in each of these infinitely many
parts possesses *as many* (*tosa*) parts *as (hosa) of* any other part, namely *the same finite number*
of parts *(Parmenides* 144d4-e3),
exactly as with the self-similar One and the Finite in Zeno's Fragment B3,
(d2) is a paradigmatical Platonic intelligible Being, and
(d3) coincides with Zeno's true Being.
Thus *Zeno's true Being is also the philosophic analogue of periodic anthyphairesis*,
something that can only mean that
Zeno's true Being is a philosophic imitation of the Pythagorean discovery of the infinite periodic
anthyphairesis of the diameter to the side of a square. Since the only way that Zeno would have
as principles for his true Beings the anthyphairetic Infinite and the Finite is by imitating the
Pythagoreans, it follows the Pythagorean proof of incommensurability has occurred at an early
date, namely before Zeno formed his paradoxes and arguments.
For details on our anthyphairetic interpretation of Zeno's paradoxes and arguments the reader is
referred to Negrepontis 2023.

**Section 14**. *The second confirmation of our reconstruction:*
*our interpretation of the Pythagorean definition of the three kinds of Angles at the position*
*II.13/14.a and of the Postulate 4 of the* Elements *at the point II.13/14.b, in terms of the two*
*Pythagorean principles Infinite and Finite*

In the long fascinating passage Proclus, *In Eucliden* 131.13-134.7 describes how the Pythagoreans maintained that the cause of the right angle is the Pythagorean principle of the Finite, and the cause of the acute and obtuse angles is the Pythagorean principle of the Infinite. We will outline our interpretation of the Pythagorean approach to the three kinds of angles, not by simply defining them, but by appealing to the two Pythagorean principles.

**14.1**. *The cause of the obtuse and acute angles*
*on the one hand is said by Proclus to be the Pythagorean principle of the Infinite, and*
*on the other hand by our interpretation is*
*the infinite anthyphairesis of the diameter to the side of a square,*
*suggesting that*
*the Pythagorean principle of the Infinite*
*coincides with*
*the infinite anthyphairesis of the diameter to the side of a square*

We will argue that
Propositions II.12, 13 for isosceles triangles express the inequalities characterzing the acute and obtuse angles. and
the isosceles Pythagorean theorem II.4/5 expresses the equality characterizing the right angle.

The Proclus passage describes how the acute and obtuse angles are being caused the the Pythagorean principle of the Infinite. The principle of the Infinite appears in various expressions, such as *apeiron, aperantos, aoriston* (indefinite), *mallon-kai hetton* (more-and-less), *meizon-kai -ellason*(greater-and-smaller), *ametria* (incommensurability), infinite motion and change, *ep'apeiron proodos* (infinite progress). These expression clearly point to an anthyphairetic infinite, e.g. the more-and-less or the greater-and-smaller are the expressions used by Plato in the *Philebus* 23b-25e to exthe philosophical analogue of the anthyphairetic infinite (cf. Negrepontis, 2005).

[1] The Pythagorean principle of the Infinite has been interpreted, in Section 12, as the infinite anthyphairesis of the diameter d to the side s of the square.

[2] The infinite anthyphairesis of the diameter d to the side s of a square is equivalently described by the infinite double sequence of the side $p_n$ and diameter $q_n$ numbers, n=1,2,…..
(since by Proposition II.10/11,
Anth (d,s)   = [1,2,2,2,......], and
Anth ($q_n,p_n$) = [1,2,2,2,...,2 (n-1 times)] for n=1,2,...).

[3] The sequence of the side and diameter numbers satisfies
the Pell property $q_n^2=2p_n^2+(-1)^n$ (Proposition II.10/11.b).

[4] The angle $\omega_n$, given by

*Definition.* We let $\omega_n$ be the angle contained by the equal sides $p_n$, $p_n$, of the isosceles triangle with sides $p_n$, $p_n$, $q_n$.

is determined by the side and diameter numbers, an thus by the anyphairesis of the diameter to the side of a square. The Pell property of the side and diameter numbers, together with Propositions II.12,13 applied to the isosceles triangle with sides $p_n, p_n, q_n$, show the following

*Proposition.* The angle $\omega_{2n-1}$ is acute, and the angle $\omega_{2n}$ is obtuse for every n=1,2,... .

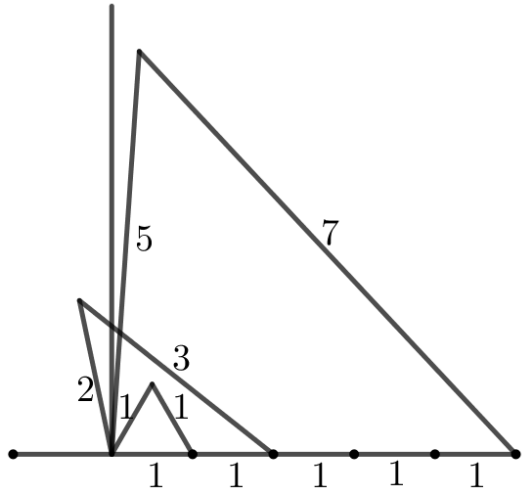

*Proof.*
By *the Pell property,*
$q_{2n-1}^2 = 2 p_{2n-1}^2 - 1$, hence $q_{2n-1}^2 < 2 p_{2n-1}^2$,
hence, by Propositions II.4/5, II.12, and II.13, the angle $\omega_{2n-1}$ contained by the equal sides $p_{2n-1}$, $p_{2n-1}$ of the isosceles triangle with sides $p_{2n-1}$, $p_{2n-1}$, $q_{2n-1}$ is *acute*; and

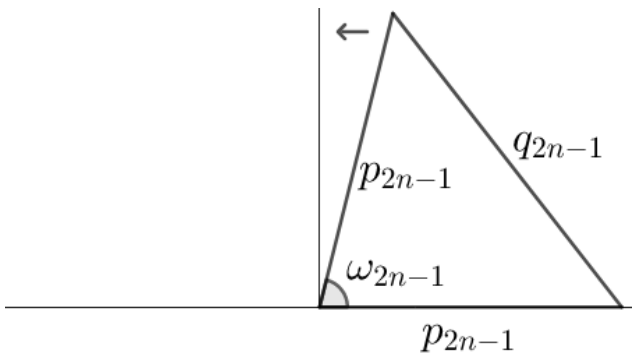

By *the Pell property,*
$q_{2n}^2 = 2p_{2n}^2 + 1$, hence, $q_{2n}^2 > 2p_{2n}^2$,
hence, by Propositions II.4/5, II.12, and II.13, the angle $\omega_{2n}$ contained by the equal sides in the isosceles triangle with sides $p_{2n}$, $p_{2n}$, $q_{2n}$ is *obtuse*.

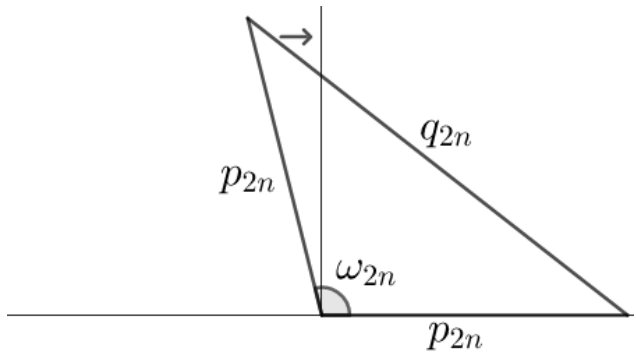

[5] Both the increasing sequence of acute angles $\omega_{2n-1}$ and the decreasing sequence of obtuse angles $\omega_{2n}$ converge to the right angle as n goes to infinity.

| $q_1,p_1$ | | | | | | | | | |
|---|---|---|---|---|---|---|---|---|---|
| | | | | | | | | | $q_2,p_2$ |
| | $q_3,p_3$ | | | | | | | | |
| | | | | | | | | $q_4,p_4$ | |
| | | $q_5,p_5$ | | | | | | | |
| | | | | | | | $q_6,p_6$ | | |
| | | | $q_7,p_7$ | | | | | | |
| | | | | | | $q_8,p_8$ | | | |
| | | | | … | | | | | |
| | | | | | … | | | | |
| | | | | $d^2-2s^2$ | | | | | |

Indeed, the diameter $q_n$ and the side $p_n$ of the isosceles triangle tends to infinity, and informally we have that
 $q_n^2/p_n^2 = 2/1 - (1/p_n^2)$,
and thus the sequence of angle ($\omega_n$) approaches and converges to the right angle. Hence, the terms of the sequence ($\omega_{2n-1}$) of acute angles eventually exceed every acute angle, and the terms of the sequence ($\omega_{2n}$) of obtuse angles eventually become smaller than every obtuse angle.

[6] This allows for a definition of the two kinds of angles, the obtuse and the acute angle, ultimately in terms of the infinite anthyphairesis of the diameter to the side of a square, as follows:

*Pythagorean definition of the acute and the obtuse angle.*
 An angle $\omega$ is *acute* if there is an *odd* number 2n-1, such that $\omega < \omega_{2n-1}$, and
 an angle $\omega$ is *obtuse* if there is an *even* number 2n, such that $\omega > \omega_{2n}$.

*Note*. At this point it is explained why the two Propositions II.12,13, although proved by means of the Pythagorean theorem and the Gnomonic Propositions II.4 and II.7, appear in Book II *later than expected*, in fact after the two Propositions II.9 and 10, intimately connected with the Pell property of the side and diameter numbers.
The reason is that Propositions II.12 and 13 are first needed, together with *the Pell property* for the side and diameter numbers (Proposition II.10/11.c), to derive the fact that the Pythagoreans were producing the two of the three kinds of angle, the obtuse and the acute angle, from the Pythagorean principle of the Infinite.
The position of Propositions II.12 and 13 suggests that we apply these two propositions into the side and diameter numbers in order to obtain acute and obtuse angles.

[7] Thus Proclus' claim that
*the Pythagorean principle of the Infinite* is the cause of the acute and the obtuse angles
is confirmed if we interpret the Pythagorean principle of the Infinite as
the infinite anthyphairesis of the diameter to the side of a square.

[8] And this in turn confirms our interpretation, in Section 6.12, that *the Pythagorean principle of the Infinite* coincides with the infinite anthyphairesis of the diameter to the side of a square.

**14.2**. *The cause of the right angle on the one hand is said by Proclus to be the Pythagorean principle Finite, and on the other hand by our interpretation is the preservation of the Gnomons in the anthyphairesis of the diameter to the side of a square, suggesting that the Pythagorean Finite coincides with the preservation of the Gnomons in the anthyphairesis of the diameter to the side of a square,*

Proclus claims that

*the right angle is produced from the principle of the Finite.*

The principleof the Finite is described by the terms peras (finite), horos, horismos (limit).

[1] The Pythagorean principle of the Finite has been interpreted, in Section 12, as
the preservation of the Gnomons in the anthyphairesis of the diameter to the side of a square; in the notation of Section 11

$$b^2=2bc_1+c_1^2=G_1,$$
$$c_1^2=2c_1c_2+c_2^2=G_2.$$

[2] By Proposition 3.2,
$bc_2=c_1^2$,
and since $c_1=a-b$, $c_2=b-2c_1=b-2(a-b)=3b-2a$, we have, by II.7,
$b(3b-2a)=(a-b)^2=a^2+b^2-2ab$,
$3b^2-2ab= a^2+b^2-2ab$,
**$2b^2=a^2$.**
Thus the preservation of the Gnomons results in the equality $a^2=2b^2$, which,
by the Pythagorean theorem II.4/5 and its converse II.12 & 13,
is equivalent to the right angle.

[3] Thus Proclus' claim that
*the Pythagorean principle of the Finite* is the cause of the right angle
is confirmed if we interpret the Pythagorean principle of the Finite as
the preservation of the Gnomons in the anthyphairesis of the diameter to the side of a square.

[4] And this in turn confirms our interpretation, in Section 12, that
*the Pythagorean principle of the Finite*
coincides with the preservation of the Gnomons in the anthyphairesis of the diameter to the side of a square.

**14.3.** *Note*. Morrow, 1970, p. 104, writes in his Note 93 to Proclus' 129.17:

> A cryptic reference to the distinction between right angles and obtuse or acute angles. See 131.13-134.7

This remark shows that he has not understood the connection of the Pythagorean treatment of the three kinds of angles with the anthyphairetic proof of incommensurability of diameter to the side of a square, and the resulting side and diameter numbers and their Pell property. To the best of our knowledge our interpretation is new.

**14.4**. The meaning that the Pythagoreans assigned to
*the Fourth Postulate* in Euclid's *Elements: All right angles are equal*
*according to* Proclus, *In Eucliden* 191,5-15

There is something anomalous with the status of Postulate 4. Proclus, *In Eucliden* 188,11-189,12, gives a quite simple proof of Postulate 4; thus, in Euclid's scheme, Postulate 4 is not really a Postulate, but a Proposition with a proof. Proclus' passage describing the establishment of Posstulate 4 from the Pythagorean and anthyphairetic principle of the Finite, strikes us as a Pythagorean form of Thales theorem. We describe our interpretation in Section 14.6, below.

**14.5**. Proclus, *In Eucliden* 191,5-15

[1] Φανερὸν δὲ καὶ ἐκ τοῦδε *τοῦ αἰτήματος*, ὅτι
*ἡ ὀρθότης* τῆς γωνίας
τῇ ἰσότητι
*συγγενής* ἐστιν,
ὥσπερ *ἡ ὀξύτης καὶ ἀμβλύτης*
τῇ ἀνισότητι. 191,5-7

*This postulate* also shows that
rightness of angles
is akin to equality,
as acuteness and obtuseness
are akin to inequality.

[2] καὶ γάρ ἐστιν
*ἡ μὲν ὀρθότης* αὐτῇ
τῇ ἰσότητι
*σύστοιχος*
—ἀμφότεραι γὰρ ὑπὸ τὸ πέρας,
ὥσπερ δὴ καὶ *ἡ ὁμοιότης*—
*ἡ δὲ ὀξύτης καὶ ἀμβλύτης*
τῇ ἀνισότητι,

καθάπερ καὶ *ἡ ἀνομοιότης·*
*ἀπειρίας γὰρ ἔκγονοι* πᾶσαι.

διὸ καὶ
οἱ μὲν τὸ ποσὸν ὁρῶντες τῶν γωνιῶν
*τὴν ὀρθὴν ἴσην τῇ ὀρθῇ* λέγουσιν,
οἱ δὲ τὸ ποιὸν *ὁμοίαν.*
ὅπερ γάρ ἐστιν ἐν ποσοῖς *ἡ ἰσότης*
τοῦτο ἐν τοῖς ποιοῖς *ἡ ὁμοιότης*.
191,7-15

Rightness is in the same column with equality,
for both of them belong under the Finite,
as does likeness also.

But acuteness and obtuseness
are akin to inequality,
as is also unlikeness;
for all of them are the offspring of the Unlimited,

This is why
those who look at the quantity of angles
say that a right angle is equal to a right angle,
while others, looking at the quality, say it is similar.
For similarity has the same position among qualities
that equality has among quantities.
Proclus, *In Eucliden* 191,5-15

**14.6**. *Our interpretation of Proclus' passage on Postulate 4*

*The first difference* of the right angle from the class of acute and abtuse angles is the difference between *equality and inequality,* and this difference is a consequence of the fact that cause of the the right angle is the principle of Finite, while the case of the acute and abtuse angles is the principle of the Infinite (in [1] 131,9-17, [2] 131,17-132,6, [3] 132,6-17, [10] 191,5-7, [11]).

*Postulate 4*. All the right angles are equal

*Proof* (from the Pythagorean principle of the Finite)
Let $\omega_1$, $\omega_2$ be two right angles.
We construct two isosceles right triangles, one with right angle $\omega_1$ and with hypotenuse a and sides b, the other with right angle $\omega_2$ and with hypotenuse c and sides d.
For example, assume that $a^2=2b^2$, and set c=a+2b, d=a+b. as in the Elegant theorem.

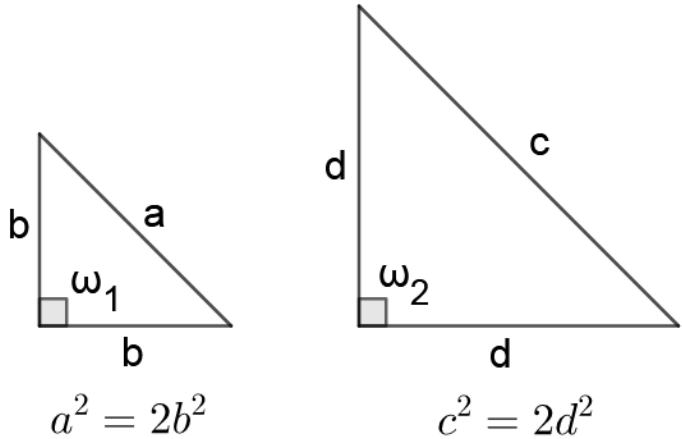

$$a^2 = 2b^2 \qquad c^2 = 2d^2$$

By the Pythagorean theorem II.4/5, $a^2=2b^2$, $c^2=2d^2$.
Then by Proposition 11.2, ad=bc.
By Proposition 9.4.1, Anth(a,b)=Anth(c,d)=[1, period(2)].
We thus see that indeed the equality of the infinite anthyphairesis follows from the Pythagorean principle of the Finite, namely the preservation of Gnomons/Application of Areas in Excess from the second to the third step of the anthyphairesis of the diameter to the side of a square.
Thus, by *equality of all right angles (Postulate 4)* the Pythagoreans mean that in *every isosceles right triangle* the anthyphairesis of the hypotenuse to the side is *equal to* [1, period(2)].
The Postulate is thus seen as a form of Thales Theorem without ratios of magnitudes.

**Section 15**. *Our answers to the Questions about Pythagorean Mathematics and Philosophy*

At this point we feel justified to declare that our reconstruction of the original Pythagorean proof of incommensurability, given in Section 11, should no longer be regarded as tentative, but reasonably firmly established.

**-Q1**. *What is the method of the original proof of incommensurability by the Pythagoreans?*

*Answer to Q1*. It is seen that Book II is divided into two parts, one before Pythagorean incommensurability and the other after, and suggests an exact position for the Pythagorean Proposition of proof of incommensurability, at point II.8/9.a, right in the middle of these two parts, together with the tools that are to be used. This suggestion is similar to the earlier suggestion by which we were led to the reconstruction of the original Pythagorean proof of the Pythagorean Theorem at Proposition II.4/5, by Plato's *Meno* 84d3-85b7. Proclus' passage on the three kinds of angles is of great importance and usefulness. It establishes that the Pythagorean proof of incommensurability is anthyphairetic and uses Application of Areas and Geometric Algebra.

**Q2**. How did the Pythagoreans conceive of the side and diameter numbers?

*Answer to question Q2*. The Pythagoreans defined them first as the unique relatively prime numbers, whose anthyphairesis is equal to a finite initial segment of the infinite anthyphairesis of the diameter to the side of a square.

**Q3**. What is the role of Proposition II.8 in Book II? Where is it used?

*Answer to Q3. Proposition II.8* is used in our reconstruction of the Pythagorean incommensurability proof.

**Q4**. Why is the converse of the Pythagorean Theorem so late in Book II?

*Answer to Q4.* The converse of the Pythagorean theorem is delayed because it is proved, exactly where it is needed and used, namely for the proof of II.13/14.a, the Pythagorean approach to the three kinds of angle.

**Q5**. Do Propositions II.5/6, 6/7 really belong in Book II? Do Propositions II,5/6, 6/7 indicate that the Pythagoreans were practicing Geometric Algebra? And if so, is there a Pythagorean use of Application of Areas II.5/6 or II.6/7?

*Answer to Q5. Geometric Algebra,* and *Application of Areas* are used in an algebraic manner in the original Pythagorean proof of incommensurability, specifically for the computation of the second and third anthyphairetic quotients. This appears to be the only case of use of the Application of Areas in Excess by the Pythagoreans. Our reconstruction of the original Pythagorean proof is confirmed (a) by our anthyphairetic interpretation of the two Pythagorean principles, Infinite and Finite/Finitizing, and (b) by the definition of the three kinds of angles and the derivation of Postulate 4 by the Pythagoreans, in terms of the Pythagorean principles Infinite and Finite.
Raising the tentative nature of our reconstruction of the Pythagorean incommensurability implies that we raise as well the tentative nature of the Pythagorean use of Geometric Algebra.

**Q6**. Did the Pythagoreans know proof by mathematical induction, and did they prove the Pell property of the side and diameter numbers, or did they simply verify a few initial cases of the Pell property?

*Answer to Q6.* We have argued earlier on the basis of the accounts by Theon, Iamblichus, but mainly by Proclus' *Commentary to the Republic,* where we have noted two striking linguistic similarities that suggest an inductive proof of the Pell property by means of the Proposition II.10. Proposition II.13/14.a provides a strong independent argument that the Pythagoreans had indeed a proof of the Pell property of the side and diameter numbers, since now the Pythagoreans actually *use* the Pell property in order to define the three kinds of angles from the higher principles of Infinite and Finite. *It is safe to assume that the Pythagoreans would not use the Pell property, if they had not proved it rigorously.* Hence the very presence of this passage shows that, despite the claims of Mueller, 1981, Unguru, 1991, Fowler, 1991, the Pythagoreans not only had knowledge of induction, but in fact employed induction in order to prove the Pell property of the side and diameter numbers. In fact we may regard symbolically the Pythagorean proof of the Pell property as *the birth of mathematical induction*.

**Q7**. What is the meaning of the two Pythagorean principles Infinite and Finite?

*Answer to Q7.* Propositions II.13/14.a,b confirm clearly that the two Pythagorean principles, the Infinite and the Finite are to be interpreted as Infinite anthyphairetic division, and "preservation" of Gnomons/Application of Areas, respectively. That Pythagorean infinite is anthyphairetic can also be seen from the Philolaus *Fragment* 6, where the Infinite can hardly mean anything else besides the Infinity of the musical anthyphairesis (Section 2.5).

**Q8**. *Dating of the discovery of the Pythagorean incommensurability*

*Answer to Q8.* Both the association of the Pythagorean discovery with Hippasus and more convincingly our anthyphairetic interpretation of Zeno's arguments and paradoxes, in particular of Zeno's *Fragment* B3, show that Zeno's philosophy was shaped in imitation of the Pythagorean anthyphairetic proof of the incommensurability of the diameter to the side of a square, and thus the Pythagorean discovery is prior to Zeno's time.

**Q9**. *Why Euclid did not retain the original Pythagorean proof .of incommensurability in Book II of his* Elements*?*

*Answer to Q9.* It should be clear that Euclid made many changes in Book II. Thus the Pythagorean Theorem was moved from II.4/5 to I.47, the converse Pythagorean from II.12 & 13 to I.48, the Application of Areas in Defect II.5/6 and in Excess II.6/7, to the more general Propositions VI.28 and VI.29, respectively, but retaining II.11, because it was needed before Book VI for the construction in Book IV of the canonical pentagon. In line withese changes he omitted the Pythagorean proof of incommensurability of the diameter to the side of a square, which we have found to in position II.8/9, because of the presence of a more general theorem in the problematic Proposition X.9.

As a final conclusion, our approach unifies the restored/completed Book II of the *Elements,* which in its Euclidean form appears as a collection of mostly unrelated Propositions: it emerges that the principal original purpose of Book II is *to develop the tools and provide the proof for the incommensurability of the diameter to the side of a square* (and also for the incommensurability of the mean and extreme ratio following from Proposition II.11), and *to provide the basic properties of the side and diameter numbers*.

*Transalations*

Iamblichus, *On the Pythagorean Life, translated with Notes and Introduction by G. Clark*, Liverpool University Press, 1989.
Jamblique, *In Nicomachi arithmeticam, Edition and translation by N. Vinel*, Fabrizio Serra Editore, Pisa and Rome, 2014.
Pappus, *Commentary on Book X of Euclid's Elements, Arabic text and translation by W. Thomson; Introductory Remarks, Nores, and a Glossary of technical terms by G. Junge and W. Thomson,* Cambridge, 1930
Plato, *Republic, translated by P. Shorey.* Harvard University Press, Cambridge, Massachusetts; W. Heinemann, London, 1969.
Plato, *Parmenides, Theaetetus, Sophist, Politicus, Philebus, Epinomis*, *translated by H. N. Fowler.* Harvard University Press, Cambridge, Massachusetts; W. Heinemann, London, 1921.
Proclus, *A Commentary on the First Book of Euclid's Elements. translated, with Introduction and Notes* by G. R. Morrow, Princeton University Press, Princeton, 1970.
Proclus, *Commentary on Plato's Republic*, Volume II, Essays 7-15, translated with an introduction and notes by D. Baltzly, J.F. Finamore, G, Miles, Cambridge University Press, 2022
Plutarch, *Lycurgus and Numa, translated by B. Perrin*, William Heinemann Harvard University Press Cambridge, Massachusetss, 1967.
Theon of Smyrna, *Mathematics Useful for Understanding Plato*, *Translated from the 1892 Greek/French edition of J. Dupuis by R. and D. Lawlor*, Wizards Bookself, San Diego. 1979.

*Note.* Unless otherwise indicated, the translations of the ancient Greek passages in the paper follow the translations given in this list, with modifications by the authors.

Authors address:
Stelios Negrepontis snegrep@snegrep@gmail.com
Vassiliki Farmaki vfarmaki@math.uoa.gr

Department of Mathematics, Athens University, Athens, Greece
January 15, 2025